\documentclass[journal,final,twocolumn]{IEEEtran}
\usepackage{graphicx}
\usepackage{subfigure}
\usepackage[colorlinks]{hyperref}
\hypersetup{
    linkcolor=blue,
    filecolor=magenta,      
    urlcolor=teal,
    citecolor=magenta
    }
\usepackage{comment}
\usepackage{breakurl}
\usepackage{amssymb}
\usepackage{mathtools}
\usepackage{algorithm}
\usepackage{algorithmic}

\usepackage{amsfonts}
\usepackage{amsmath}
 
\usepackage{amsthm}
\theoremstyle{definition}

\usepackage{enumitem}
\usepackage{bm} 
\usepackage{algorithm}
\usepackage{algorithmic}
\usepackage[]{vaucanson-g}

\begin{document}

\title{Model-based Prediction and Optimal Control of Pandemics by Non-pharmaceutical Interventions}
%////////////////////////////////////////////
\author{Reza~Sameni$^\textnormal{*}$,~\IEEEmembership{Senior~Member,~IEEE}%
\thanks{Manuscript received May 30, 2021; revised October 19, 2021; accepted November 15, 2021.}%
\thanks{\textsuperscript{*}R.~Sameni is with the Department of Biomedical Informatics, Emory University, Atlanta, GA (e-mail: \url{rsameni@dbmi.emory.edu}).}%
\thanks{Copyright (c) 2021 IEEE. Personal use of this material is permitted.  However, permission to use this material for any other purposes must be obtained from the IEEE by sending an email to pubs-permissions@ieee.org.}%
}%\newline
%////////////////////////////////////////////

\markboth{Published in IEEE Journal of Selected Topics in Signal Processing: \url{https://doi.org/10.1109/JSTSP.2021.3129118}}{R. Sameni: Model-based Prediction and Control of Pandemics}
\maketitle
%////////////////////////////////////////////
\begin{abstract}
A model-based signal processing framework is proposed for pandemic trend forecasting and control, by using non-pharmaceutical interventions (NPI) at regional and country levels worldwide. The control objective is to prescribe quantifiable NPI strategies at different levels of stringency, which balance between human factors (such as new cases and death rates) and cost of intervention per region/country. Due to infrastructural disparities and differences in priorities of regions and countries, strategists are given the flexibility to weight between different NPIs and to select the desired balance between the human factor and overall NPI cost.

The proposed framework is based on a \textit{finite-horizon optimal control} (FHOC) formulation of the bi-objective problem and the FHOC is numerically solved by using an ad hoc \textit{extended Kalman filtering/smoothing} framework for optimal NPI estimation and pandemic trend forecasting. The algorithm enables strategists to select the desired balance between the human factor and NPI cost with a set of weights and parameters. The parameters of the model are partially selected by epidemiological facts from COVID-19 studies, and partially trained by using machine learning techniques. The developed algorithm is applied on ground truth data from the Oxford COVID-19 Government Response Tracker project, which has categorized and quantified the regional responses to the pandemic for more than 300 countries and regions worldwide, since January 2020. The dataset was used for NPI-based prediction and prescription during the XPRIZE Pandemic Response Challenge.
\end{abstract} % 
%////////////////////////////////////////////
%////////////////////////////////////////////
\IEEEpeerreviewmaketitle
%////////////////////////////////////////////
\section{Introduction
\label{sec:introduction}}
The COVID-19 pandemic highlighted the fact that social life--- as a \textit{dynamic system}--- is always in a \textit{metastable} condition, which is continuously prone to pandemic outbreaks (regardless of the severity or geographical origin of pandemics). Parallel to medical solutions and vaccinations against known viruses, the rapid and effective response to future pandemics requires proactive planning and interdisciplinary collaborations between different scientific communities. Specifically, some of the prominent contributions, which can be made by the signal processing and data science communities include: 1) developing accurate spatio-temporal forecasting models at different levels of abstraction, which simulate pandemic outbreaks and trends, 2) identifying quantifiable \textit{non-pharmaceutical intervention} (NPI) plans, with fact-based estimates of the impact and cost of each NPI \cite{Hale2021}, and 3) simulated multi-objective pandemic response strategies, which balance between NPI cost and effectiveness, to help governments and policymakers in resource allocation and fact-based decision making to control new pandemic waves.

In this context, NPIs refer to actions and policies adopted by individuals, authorities or governments that help slowing down the spread of epidemic diseases. Enforcement of social distancing, face covering, restrictions on social events and public transportation, etc. were among the NPIs that were experienced during the COVID-19 pandemic. NPIs are among the best ways of controlling pandemic diseases when vaccines or medications are not yet available\footnote{See Centers for Disease Control and Prevention guidelines on NPIs: \url{ https://www.cdc.gov/nonpharmaceutical-interventions/}.}.

\begin{figure*}
    \centering
    \includegraphics[angle = 90, trim=2.4in 1.5in 2.2in 0.5in,clip,width=1.4\columnwidth]{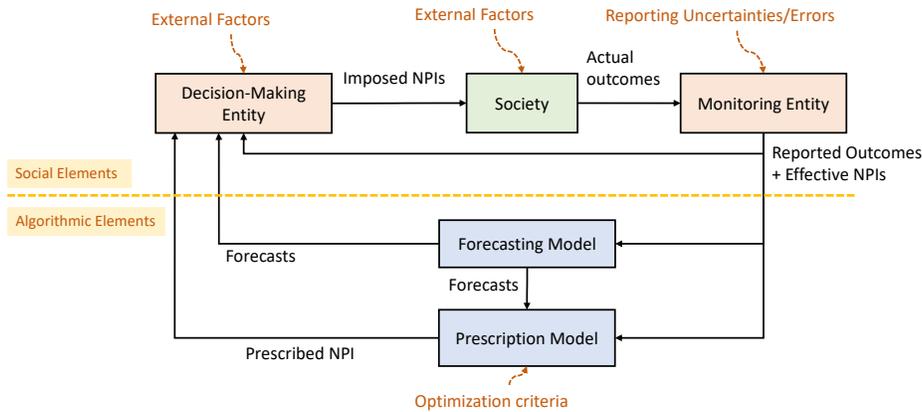}
    \caption{A closed-loop system representation for pandemic control via non-pharmaceutical interventions (NPIs). Social constraints imposed by governments, states, etc. is influenced by the pandemic trend and external political and economic factors. Beyond the level of social interactions, the pandemic spread is also influenced by external factors such as vaccination, new virus variants, etc. The NPIs \textit{imposed} by decision-makers practically differ from the actual \textit{effective} NPIs practiced by the society. The monitoring entities (healthcare system and other organizations) observe the ``ground truth'' effective NPIs, the number of new cases, hospitalized and death cases, which contain inevitable errors and uncertainties. Data analysts contribute to this loop by forecasting the future trend of the pandemic and providing fact-based prescriptions of future NPI to control the pandemic, while satisfying other socioeconomic objectives.}
    \label{fig:NPISystem}
\end{figure*}

During the COVID-19 pandemic, several attempts were made to categorize and quantify the various NPIs of different regions and nations, which was essential for comparing the effectiveness of regional policies in containing the pandemic spread. By using signal processing and machine learning techniques, the quantified NPI can be used to forecast the future trends of the pandemic and to simulate \textit{``what if scenarios''} for the better management of human and medical resources, and to eventually prescribe appropriate NPI for controlling the pandemic \cite{miikkulainen2021prediction}. The overall pandemic monitoring, forecasting and control cycle can be seen as a closed-loop system consisting of social and algorithmic elements, as shown in Fig.~\ref{fig:NPISystem}.

The Oxford COVID-19 Government Response Tracker (OxCGRT) was one of the NPI tracking projects, which were launched and regularly updated during the COVID-19 pandemic \cite{OxCGRT2020}. OxCGRT has been used by the data science community for NPI-based prediction and prescription planning. Specifically, the XPRIZE Pandemic Response Challenge addressed the problem of predicting future trends of the pandemic in different regions and countries, and prescribing efficient NPIs that compromise between the number of new cases and a weighted-cost of intervention \cite{XPRIZEPandemicResponse2020,miikkulainen2021prediction}. The challenge was motivated by the fact that due to disparities in infrastructures, available resources and priorities, policymakers worldwide tend to give different weights to each NPI and are interested to know the impacts and consequences of each policy in advance. During this challenge, the Alphanumerics Team from the Department of Biomedical Informatics at Emory University, adopted a model-based signal processing approach, based on estimation theory and \textit{finite-horizon optimal control}, to address the problem of weighted NPI prescription. Our team was among the finalists of this challenge.

The notion of multi-objective finite-horizon pandemic control has also been considered by other researchers in simulated scenarios \cite{Mallela2020,GUAN2020394}, and also on real data \cite{Perkins2020}. The latter is closely related to the hereby developed bi-objective optimization via NPI. Over the past year, several other control-theoretical approaches have emerged for pandemic control \cite{Djidjou-Demasse2020, Patterson-Lomba2020, Shah2020}, which have been mainly tested in simulated scenarios or over regional data.

In order to develop a general framework, which is applicable for all regions worldwide, the OxCGRT dataset was used for the current study. Since the only globally registered NPIs in this dataset are the total confirmed cases, the total confirmed deaths and the daily NPIs, we adopted an extension of a generic \textit{susceptible-infected} (SI) compartmental model from our previous work, as the base model for all regions/countries \cite{sameni2020mathematical}. The proposed model parameters are trained on historic data and used to predict future trends from input NPI by using an extended Kalman filter. It is further shown that the forecasting model can be integrated with a finite-horizon optimal controller to find the optimal daily NPIs with arbitrary NPI cost weight vectors. It should be noted that since the XPRIZE Challenge was held before the global availability of COVID-19 vaccines, the interventions considered for all nations do not include vaccinations. This fact is also reflected in the compartmental model detailed in the sequel. Nevertheless, the proposed framework is generic and can be extended to other compartmental models, when more accurate data exist for a specific region/nation. 

While the majority of recent pandemic data analysis research--- including the winners of the XPRIZE Pandemic Response Challenge \cite{Lozano2021,Janko2021}--- are based on data-driven machine learning (ML) techniques and deep neural networks, this research demonstrates how classical signal processing and optimal control theories can be used and combined with ML techniques for solving pandemic trend forecasting and NPI prescription. The advantages of the proposed approach are manifold, including: theoretical support for the predictions/prescriptions, ease of interpretation through fact-based data models, quantitative performance bounds, and lower computational/training cost as compared with fully data-driven ML-based approaches.

In Section \ref{sec:dataset}, the OxCGRT NPI database is explained. Section \ref{sec:datamodel} details the data model. The developed finite-horizon optimal NPI prescription framework is elaborated in Section \ref{sec:finitehorizoncontrol}. This framework is combined with an extended Kalman filter/smoother for pandemic forecasting in Section \ref{sec:prediction_prescription}, followed by the details of model training and implementation in Section \ref{sec:modeltraining}. The results on real data from the COVID-19 pandemic and a detailed discussion on the proposed method are presented in Sections \ref{sec:results} and \ref{sec:discussion}, followed by concluding remarks and future directions. The source codes of the developed models and algorithms are provided online for reference \cite{EpidemicModelingCodesSameni2020}.

\section{The non-pharmaceutical interventions dataset}
\label{sec:dataset}
To date, the Oxford COVID-19 Government Response Tracker (OxCGRT) is an ongoing project \cite{Hale2021}, which categorizes and quantifies the NPI policies of different regions/nations since the beginning of the pandemic. The dataset was used in the XPRIZE Pandemic Response Challenge \cite{XPRIZEPandemicResponse2020}, which addressed the problem of pandemic trend forecasting in different regions/countries under social interventions (e.g., social distancing, mandatory mask wearing, social gathering prohibitions, closure of schools and public transportation limitations), and prescribing efficient NPIs that compromise between human factors (infection and death rates) and a weighted cost of intervention. The subset of OxCGRT NPI categories used by the XPRIZE Challenge and the current study are listed in Table~\ref{tab:OxCGRTNPI}. Note that the OxCGRT dataset uses the Johns Hopkins Coronavirus dataset for US states \cite{Dong2020}, which also provides state-level information for the US. 
In OxCGRT, the level of NPI stringency for different social activities are quantified by integer values ranging from zero (minimum or no stringency), to two or four (for a maximum stringency, depending on the index). The total social stringency can be considered as a weighted combination of the different indexes. Throughout the XPRIZE Challenge, the weight vector was considered as a design parameter, which enables strategists to prioritize certain NPIs over the others.

\begin{table*}[tb]
    \centering
\caption{Subset of NPI indexes from the Oxford COVID-19 Government Response Tracker project used in this study, adopted from \cite{Hale2021}}
    \begin{tabular}{|c|l|p{12cm}|}
    \hline
      Index & Description & Values\\\hline
       C1 & School closing & 0: no measures 
1: recommend closing or all schools open with alterations resulting in significant differences compared to non-COVID-19 operations, 2: require closing (only some levels or categories, e.g., only high school or public schools), 3: require closing all levels, Blank: no data\\\hline
       C2 & Workplace closing & 0: no measures, 1: recommend closing (or recommend work from home), 2: require closing (or work from home) for some sectors or categories of workers, 3: require closing (or work from home) for all-but-essential workplaces (e.g., grocery stores, doctors), Blank: no data\\\hline
       C3 & Cancel public events & 0: no measures, 1: recommend canceling, 2: require canceling, Blank: no data\\\hline
       C4 & Restrictions on gatherings & 0: no restrictions, 1: restrictions on very large gatherings (above 1000 people), 2: restrictions on gatherings (101--1000 people), 3: restrictions on gatherings (11--100 people), 4: restrictions on gatherings (up to 10 people), Blank: no data\\\hline
       C5 & Close public transport & 0: no measures, 1: recommend closing (or significantly reduce volume/route/means of transport available), 2: require closing (or prohibit most citizens from using it), Blank: no data\\\hline
       C6 & Stay at home requirements & 0: no measures, 1: recommend not leaving house, 2: require not leaving house with exceptions for daily exercise, grocery shopping, and `essential' travels, 3: require not leaving house with minimal exceptions (e.g., allowed to leave once a week, or only one person at a time), Blank: no data\\\hline
       C7 & Internal movement restrictions & 0: no measures, 1: recommend not to travel between regions/cities, 2: internal movement restrictions in place, Blank: no data\\\hline
       C8 & International travel controls & 0: no restrictions, 1: screening arrivals, 2: quarantine arrivals from some or all regions, 3: ban arrivals from some regions, 4: ban on all regions or total border closure, Blank: no data\\\hline
       H1 & Public information campaigns & 0: no COVID-19 public information campaign, 1: public officials urging caution about COVID-19, 2: coordinated public information campaign (e.g., across traditional and social media), Blank: no data\\\hline
       H2 & Testing policy & 0: no testing policy, 1: only those who have symptoms AND meet specific criteria (e.g., key workers, admitted to hospital, came into contact with a known case, returned from overseas), 2: testing all symptomatic people, 3: open public testing (e.g., ``drive through'' testing available to asymptomatic people), Blank: no data\\\hline
       H3 & Contact tracing & 0: no contact tracing, 1: limited contact tracing (not done for all cases), 2: comprehensive contact tracing (done for all identified cases)\\\hline
       H6 & Facial coverings & 0: no policy, 1: recommended, 2: required in some specified shared/public spaces outside the home with other people present, or some situations when social distancing not possible, 3: required in all shared/public spaces outside the home with other people present or all situations when social distancing not possible, 4: required outside the home at all times regardless of location or presence of other people\\\hline
    \end{tabular}
    \label{tab:OxCGRTNPI}
\end{table*}

\section{Data Model}
\label{sec:datamodel}
\begin{figure}[t]
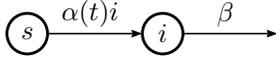

     \centering
     \MediumPicture\VCDraw{\begin{VCPicture}{(0,-1)(5,0)}
 % states
 \ShowState
 \State[s]{(0,0)}{A} \State[i]{(3,0)}{B} \HideState
 \State[R]{(6,0)}{C}
 % initial--final
 %\Initial{A} \Final{B}
 % transitions
 \EdgeL{B}{C}{\beta}
 \EdgeL{A}{B}{\alpha(t) i}
 \end{VCPicture}}
     \caption{The base susceptible-infected compartmental model with NPI-controlled infection rate}
     \label{fig:si_model}
 \end{figure}

The two major classes of methods for epidemic disease spread modeling are: 
\begin{enumerate}[leftmargin=*]
    \item \textit{Compartmental models}, which split the total population of a region into various compartments (groups) such as susceptibles, exposed, infected, recovered, vaccinated, deceased, etc. These compartments are used to form parametric differential/difference equations, which are fit on real data and are analytically or numerically solved to predict future trends of the disease spread. 
    \item \textit{Agent-based models}, which model the behaviors of individuals and their interactions at a simplified level-of-abstraction. Using these models, large groups of agents are generated in stochastic simulated environments, as they randomly move, interact and probabilistically pass the infection to one another, recover, pass away, etc. The population-level properties are obtained by ensemble averaging over the entire population.
\end{enumerate}
Each approach has its advantages and limitations. For large population sizes at regional or national levels--- which is the scope of the current study--- the first approach is asymptotically accurate and is more advantageous as it can be analytically studied in a rigorous mathematical framework and combined with state estimation techniques for forecasting, and optimal control theories for NPI prescription. Therefore, the first approach was adopted for this study, by using a contact-controlled time-variant version of the so-called \textit{susceptible-infected (SI)} compartmental model shown in Fig.~\ref{fig:si_model}. This model is a simplified variant of the general multi-compartment models studied in our previous research \cite{sameni2020mathematical}. Apparently, more accurate models can be used for the regions that additional data such as the number of recovered, hospitalized, vaccinated, or the age pyramid of the population are available. However, for the current study, since the global data provided in the OxCGRT dataset were the number of daily confirmed cases, total death cases, and the regional NPIs, the same SI model is used for all regions and countries.

The nonlinear dynamic equations corresponding to the proposed compartmental model are:
\begin{equation}
\begin{array}{l}
\dot{s}(t) = -\alpha(t)s(t)i(t)\\
\dot{i}(t) = \alpha(t)s(t)i(t) - \beta i(t)\\
\dot{\alpha}(t) = -\gamma \alpha(t) + \gamma h[\mathbf{u}(t)]
\end{array}
\label{eq:dynamicmodel}
\end{equation}
where
\begin{itemize}[leftmargin=*]
    \item $s(t)$ is the fraction of population in a region/country that is susceptible at time $t$ (i.e., the susceptible population divided by the population size $N$);
    \item $i(t)$ is the fraction of population that is infected and contagious at $t$ (i.e., the infected contagious population divided by the regional population size $N$);
    \item $\mathbf{u}(t)\in\mathbb{R}^L$ is the NPI vector considered as an exogenous control input ($L=12$, for the list of NPIs in Table~\ref{tab:OxCGRTNPI}, used for the XPRIZE Challenge \cite{XPRIZEPandemicResponse2020}). The full description of the OxCGRT data NPI set is detailed in \cite{OxCGRT2020};
    \item $\alpha(t)$ is the time-variant contagion rate, with inverse time unit;
    \item $h[\mathbf{u}(t)]$ is a causal monotonic function of the NPI, which maps the NPI to the contagion rate;
    \item $\beta$ is the rate of elimination from the contagious group (through quarantine, recovery, or death), assumed to be constant in the simplified case;
    \item $\gamma$ is the \textit{action to effect rate} (or the inverse of the NPI lag to the inter-individual contact rate), which accounts for the delay between adopting an NPI policy and the onset of its practical effectiveness in the pandemic trend. The third equation in (\ref{eq:dynamicmodel}) is equivalent to $\alpha(t) = \gamma\exp[-\gamma (t - t_0)] * h[\mathbf{u}(t)]$ (for $t \geq t_0$), which is a smoothed version of $h[\mathbf{u}(t)]$. As a corner case, $\gamma \rightarrow \infty$ represents zero latency between action and effect, resulting in $\alpha(t) = h[\mathbf{u}(t)]$.
    % , which is a causal function of $\mathbf{u}(t)\in\mathbb{R}^L$ input impact on the  
\end{itemize}
The parameters $\beta$, $\gamma$ and the function $h[\mathbf{u}(t)]$ require learning using the observed variables, as detailed in Section \ref{sec:modeltraining}. Furthermore, in \cite{sameni2020mathematical} we showed how the infection \textit{reproduction rate} $\mathcal{R}_t$ can be calculated from $\alpha(t)$ and $\beta$. Specifically, using the eigenanalysis-based definition of the reproduction rate proposed in \cite{sameni2020mathematical}, during the pandemic outbreak, when only several percents of the population are infected and \textit{herd immunity} has not been reached, we have:
\begin{equation}
    \mathcal{R}_t \approx \exp[\Delta(\alpha(t) - \beta)]
\label{eq:ReproductionRate}    
\end{equation}
where $\Delta$ is the reproduction rate generation time-unit.

Finally, for estimation purposes, the dynamic equations in (\ref{eq:dynamicmodel}) can be related to real-world reports of the fractions of \textit{new cases}: 
\begin{equation}
    n(t) = \alpha(t)s(t)i(t) + v(t),
\label{eq:observationNewCases}
\end{equation}
or through the fraction of \textit{total confirmed cases}:
\begin{equation}
c(t) = s(t_0) - s(t) + v(t),
\label{eq:observationTotalCases}
\end{equation}
where $v(t)$ is measurement noise due to case report errors (which inevitably existed during the COVID-19 global reports in all regions/nations), and $s(t_0)$ is the initial susceptible population fraction at the beginning of the pandemic (very close or equal to 1, for an un-vaccinated initial population). 

\section{Finite-horizon optimal NPI control}
\label{sec:finitehorizoncontrol}
\subsection{Cost function and problem statement}
From (\ref{eq:dynamicmodel}), the total number of new infections over an arbitrary time window $[t_0, t_1]$ is:
\begin{equation}
    J_0(\mathbf{u}) = -\int_{t = t_0}^{t_1} \dot{s}(t)\,\mathrm{d}t = \int_{t = t_0}^{t_1} \alpha(t)s(t)i(t)\,\mathrm{d}t
\label{eq:newcasescost}
\end{equation}
and the total cost of NPIs over the same time period is
\begin{equation}
    J_1(\mathbf{u}) = \int_{t = t_0}^{t_1} \mathbf{w}(t)^T\mathbf{u}(t)\,\mathrm{d}t 
\label{eq:npicost}
\end{equation}
where $\mathbf{w}(t)$ is the NPI weight vector given as input. The motivation for the user-selected weight vector $\mathbf{w}(t)$ is that the cost of intervention is different across regions. A stereotypical example considered in the XPRIZE Challenge was that ``closing public transportation may be much costlier in London than it is in Los Angeles. Such preferences are expressed as weights associated with each intervention plan dimension, given to the prescriptor as input for each region \cite{XPRIZEPandemicResponse2020}.''

With these assumptions, the optimal NPI prescription problem can be formulated as a bi-objective optimization problem, with a total cost:
\begin{equation}
    J(\mathbf{u}) = (1 - \epsilon)J_0(\mathbf{u}) + \epsilon J_1(\mathbf{u})\,\, \text{ s.t. } \mathbf{u} \in \Gamma
\label{eq:totalnpicost}
\end{equation}
where $\epsilon\in[0, 1]$ is a free parameter that compromises between the human factor ($J_0$) and the NPI cost ($J_1$), and $\Gamma$ is the set of admissible inputs:
\begin{equation}
    \Gamma = \{\mathbf{u} | \mathbf{u}^{\min}\leq \mathbf{u}(t) \leq \mathbf{u}^{\max}, \forall t\in [t_0, t_1]\} 
\label{eq:admissibleinput}
\end{equation}
where the vectors $\mathbf{u}^{\min}$ and $\mathbf{u}^{\max}$ are (element-wise) the minimum and maximum ranges of the NPI indexes in the OxCGRT dataset (the ``values'' column in Table~\ref{tab:OxCGRTNPI}). Accordingly, $\mathbf{u}^{\min}=\mathbf{0}$ corresponding to no stringency, and $\mathbf{u}^{\max} = [3, 3, 2, 4, 2, 3, 2, 4, 2, 3, 2, 4]^T$ for the maximum stringency of each NPI index.

For a given pair of design parameters $\{\epsilon, \mathbf{w}(t)\}$, the objective is to find $\mathbf{u}^*(t)$ for all $t\in[t_0, t_1]$, such that:
\begin{equation}
    J(\mathbf{u}^*) = \min_{\Gamma}(J(\mathbf{u}))
    \label{eq:optsolution}
\end{equation}

\subsection{The Pareto optimal solution}
The problem (\ref{eq:optsolution}) can be solved by \textit{finite-horizon optimization} \cite{naidu2003optimal}. In optimal control theory, the inputs which satisfy this equation are known as \textit{Pareto optimal (efficient)}. In fact, for an arbitrary weight vector $\mathbf{w}(t)$, by sweeping $\epsilon$ over $[0, 1]$, the \textit{Pareto-optimal front} of the optimization problem is found, from which pandemic strategists can select the desired free parameter $\epsilon$ that determines the desired operation point to balance between NPI effectiveness and cost (and its corresponding optimal NPI $\mathbf{u}^*(t)$ to be adopted by the country/region).

To solve (\ref{eq:optsolution}), first the corresponding \textit{Hamiltonian} function is formed \cite[Ch. 2]{naidu2003optimal}:
\begin{equation}
\begin{array}{rl}
     \mathcal{H} = & (1 - \epsilon)\alpha(t)s(t)i(t) + \epsilon \mathbf{w}(t)^T\mathbf{u}(t) \\
     & -\lambda_1(t) \alpha(t)s(t)i(t)
     \\
     & + \lambda_2(t)[ \alpha(t)s(t)i(t) - \beta i(t)]\\
     & - \gamma\lambda_3(t)\{\alpha(t) - h[\mathbf{u}(t)]\}
\end{array}
\label{eq:generalhamiltonian}
\end{equation}
where $\lambda_1(t)$, $\lambda_2(t)$ and $\lambda_2(t)$ are known as \textit{co-states}. According to \textit{Pontryagin's minimum principle}, the co-states and the optimal solution $\mathbf{u}^*$ satisfy \cite[Ch. 6]{naidu2003optimal}:
\begin{equation}
\begin{array}{l}
\dot{\lambda}_1(t) = \displaystyle-\frac{\partial\mathcal{H}}{\partial s} = [\lambda_1(t) -\lambda_2(t)-1 + \epsilon]\alpha(t)i(t)\\
\dot{\lambda}_2(t) = \displaystyle-\frac{\partial\mathcal{H}}{\partial i} = [\lambda_1(t) -\lambda_2(t)-1 + \epsilon]\alpha(t)s(t) + \beta\lambda_2(t)\\
\dot{\lambda}_3(t) = \displaystyle-\frac{\partial\mathcal{H}}{\partial\alpha} = [\lambda_1(t) -\lambda_2(t)-1 + \epsilon]s(t)i(t) + \gamma\lambda_3(t)
\\
\mathcal{H}(\mathbf{u}^*) \leq \mathcal{H}(\mathbf{u}), \quad \forall \mathbf{u}\in \Gamma
\end{array}
\label{eq:PontryaginPrinciple}
\end{equation}
When the inputs are unconstrained, the Hamiltonian minimizer input $\mathbf{u}^*$, in the last condition of (\ref{eq:PontryaginPrinciple}), can be found by solving
\begin{equation}
\bm{\nabla}_{\mathbf{u}}\mathcal{H}(\mathbf{u}^*) = \mathbf{0},
\label{eq:hamiltonianminmsolution}
\end{equation}
where $\bm{\nabla}_\mathbf{u}\mathcal{H}$ denotes the Hamiltonian gradient with respect to the input vector $\mathbf{u}$, and the condition should hold element-wise. In this case, a sufficient condition for the existence of a solution is to have $\bm{\nabla}_{\mathbf{u}}^2\mathcal{H}(\mathbf{u}^*) \succ \mathbf{0}$ (where $\bm{\nabla}_{\mathbf{u}}^2$ denotes the Hessian operator with respect to the input vector $\mathbf{u}$ and $\succ \mathbf{0}$ denotes positive-definiteness). In the constrained-input case--- as in this problem--- where the inputs are confined to the admissible set (\ref{eq:admissibleinput}), while the global solution of (\ref{eq:hamiltonianminmsolution}) might not exist or belong to the admissible set (\ref{eq:admissibleinput}), a Hamiltonian minimizer optimal input still exists. In either case, the optimal input is found as a parametric function of the costate $\lambda_3(t)$ and the other model parameters.

The parametric optimal input found from (\ref{eq:hamiltonianminmsolution}) is next combined with (\ref{eq:PontryaginPrinciple}) and (\ref{eq:dynamicmodel}) to calculate the states, using the initial conditions and appropriate boundary conditions (also known as the \textit{transversality conditions}) on the co-states and the Hamiltonian. The desired boundary conditions, which satisfy the pandemic control problem are:
\begin{equation}
    \begin{array}{ccc}
        \lambda_1(t_1) = 0, & \lambda_2(t_1) = 0, & \lambda_3(t_1) = 0.
    \end{array}
\label{eq:transversalityconditions}    
\end{equation}
The conditions in (\ref{eq:transversalityconditions}) are the general free end-point conditions of finite-horizon optimization problems, which match the objectives of the pandemic control problem. Alternative transversality conditions that can be studied within the proposed framework are \cite[Section 2.7]{naidu2003optimal}:
\begin{enumerate}[leftmargin=*]
    \item When the end-time $t_1$ is not fixed, but we require that $i(t_1)$ reaches below $i_{\max}$ by the end of the control period (\textit{infinite-horizon} scenario). This requires the additional condition: $\mathcal{H}(t_1) = 0$.
    \item Assuming that the objective of any NPI policy over a reasonable time period $[t_0, t_1]$ (long enough to make the NPIs effective) is to bring the number of active cases down to $i(t_1) \leq i_{\max}$, where $i_{\max}$ is some target fraction of active cases (ideally zero). In this case, the second condition in (\ref{eq:transversalityconditions}) can be replaced by: $\lambda_2(t_1) [i(t_1) - i_{\max}] = 0$.
    \item We require that $i(t_1)$ drops below $i_{\max}$ any time before a maximum end time $t_f$, which requires $(t_1 - t_f)\mathcal{H}(t_1) = 0$.
\end{enumerate}

\subsection{The NPI to inter-human contact map}
\label{sec:NPIhumanInteraction}
The solution of the NPI optimization problem depends on the choice of $h[\mathbf{u}(t)]$, i.e. the NPI to inter-human contact mapping model. Intuitively, $h[\mathbf{u}(t)]$ is expected to be a monotonically decreasing function of the input NPI vector $\mathbf{u}(t)$. In other words, more strict restrictions on social contact (corresponding to the higher values in Table~\ref{tab:OxCGRTNPI}) should overall reduce the person-to-person contact rates at the population level (this was the globally accepted rationale behind the social restrictions during the COVID-19 pandemic). However, the exact shape of $h[\mathbf{u}(t)]$ generally requires learning from historic data, where the monotonic decreasing assumption acts as a constraint during learning.

Based on this intuitive assumption, we study the following two cases, which lead to closed form solutions for the optimal input as functions of the model co-states.
\subsubsection{Linear regression model}
Let us take
\begin{equation}
h[\mathbf{u}(t)] = b + \mathbf{a}^T [\mathbf{u}^{\max} - \mathbf{u}(t)]
\label{eq:linear_model}
\end{equation}
where $\mathbf{a}$ is a vector of \textit{input influence weights} and $b$ is a constant bias (intercept value). The \textit{least absolute shrinkage and selection operator} (LASSO) method falls into this category. In addition, adding the constraint $\mathbf{a} \geq \mathbf{0}$ guarantees the monotonically decreasing relationship between the NPI and $\alpha$. In other words, more stringent NPI policies have a non-increasing effect on the human interactions parameter $\alpha$ (i.e. the NPI do not have any counter-impacts on the contact rates). Inserting $h[\mathbf{u}(t)]$ in (\ref{eq:generalhamiltonian}) we find:
\begin{equation}
\bm{\nabla}_{\mathbf{u}}\mathcal{H}(\mathbf{u}) = \epsilon \mathbf{w}(t) - \gamma\lambda_3(t)\mathbf{a}
\label{eq:LinearHamiltonianGradient}    
\end{equation}
Now, since $\bm{\nabla}_{\mathbf{u}}\mathcal{H}(\mathbf{u})$ is independent of $\mathbf{u}$, depending on its sign, the Hamiltonian which is a linear function of $\mathbf{u}$, admits its minimum at one of the extreme ends of the admissible input ranges (\ref{eq:admissibleinput}). This eventually results in
\begin{equation}
u^*_k(t) = 
\left\{
\begin{array}{ll}
     u^{\min}_k: & \epsilon w_k(t) > \gamma \lambda_3(t)a_k\\
     u^{\max}_k: & \epsilon w_k(t) < \gamma \lambda_3(t)a_k \end{array}
\right.
\label{eq:optimalinputlinear}    
\end{equation}
for $k = 1, \ldots, L$.
    
\subsubsection{Quadratic regression}
In the second case, we assume
\begin{equation}
h[\mathbf{u}(t)] = b + \mathbf{a}^T [\mathbf{u}^{\max} - \mathbf{u}(t)] + \frac{1}{2}[\mathbf{u}^{\max} - \mathbf{u}(t)]^T \mathbf{S} [\mathbf{u}^{\max} - \mathbf{u}(t)]
\label{eq:quadratic_model}
\end{equation}
where $\mathbf{a} \geq \mathbf{0}$ and $\mathbf{S}\in\mathbb{R}^{L\times L}$ is a positive-definite matrix. These assumptions guarantee the monotonically decreasing relationship between the NPI and $\alpha(t)$. Multivariate constrained polynomial fitting can be used to find $\mathbf{S}$, $\mathbf{a}$ and $b$ from historic NPI and case-report data\footnote{Despite the quadratic form of (\ref{eq:quadratic_model}), since it is linear in parameters $\mathbf{S}$, $\mathbf{a}$ and $b$, multivariate constrained polynomial fitting is applicable to find the unknown parameters by using historic data, i.e. model fitting over previous NPI actions adopted by different nations/states.}. Therefore, the required constraint polynomial fitting is straightforward to implement by conventional least squares solvers. In this case, we have:
    \begin{equation}
    \bm{\nabla}_{\mathbf{u}}\mathcal{H}(\mathbf{u}) = \epsilon \mathbf{w}(t) - \gamma\lambda_3(t)\{\mathbf{a} + \mathbf{S}[\mathbf{u}^{\max} - \mathbf{u}(t)]\}
    \label{eq:QuadraticHamiltonianGradient} 
    \end{equation}
    and setting $\bm{\nabla}_{\mathbf{u}}\mathcal{H}(\tilde{\mathbf{u}})=\mathbf{0}$ gives
    \begin{equation}
        \tilde{\mathbf{u}} = \mathbf{u}^{\max} - \mathbf{S}^{-1}[\frac{\epsilon\mathbf{w}(t)}{\gamma\lambda_3(t)} - \mathbf{a}]
    \end{equation}
    From the \textit{second partial derivative test}, since $\bm{\nabla}_{\mathbf{u}}^2\mathcal{H}=\gamma\lambda_3(t)\mathbf{S}$, and the fact that $\mathbf{S}$ is assumed to be positive definite, three cases may occur: 1) if $\lambda_3(t) > 0$, $\tilde{\mathbf{u}}$ is a local minimum, 2) if $\lambda_3(t) < 0$, it is a local maximum, and 3) $\lambda_3(t) = 0$ results in a saddle point. Therefore, by applying Pontryagin's minimum principle and considering the admissible input range (\ref{eq:admissibleinput}), after some algebraic simplifications we find: 
    \begin{equation}
    u^*_k(t) =
        \!\left\{
        \begin{array}{ll}
             \!u^{\min}_k: &\! \epsilon w_k(t) > \gamma \lambda_3(t)(a_k + s_k)\\
             \!\tilde{u}_k: &\! \gamma \lambda_3(t)a_k < \epsilon w_k(t) < \gamma \lambda_3(t)(a_k + s_k) \\
            \!u^{\max}_k: &\! \epsilon w_k(t) < \gamma \lambda_3(t)a_k
        \end{array}
        \right.
    \label{eq:optimalinputquadratic}       
    \end{equation}    
where $\tilde{u}_k$ and $s_k$ are the $k$\textsuperscript{th} entries of the vectors $\tilde{\mathbf{u}}$ and $\mathbf{S}(\mathbf{u}^{\max} - \mathbf{u}^{\min})$, respectively. Comparing (\ref{eq:optimalinputquadratic}) and (\ref{eq:optimalinputlinear}), it is clear how the quadratic case simplifies to the linear case when $\mathbf{S} \rightarrow \mathbf{0}$. It is also seen that in the linear case, the ``optimal NPI'' is always one of the extreme cases $u^{\min}_k$ (no action) or $u^{\max}_k$ (maximum stringency). But when $h(\cdot)$ is nonlinear, intermediate interventions may also be in the optimal NPI set.  
\section{A unified pandemic trend predictor and NPI prescriptor}
\label{sec:prediction_prescription}
For an ideal deterministic model, the state and co-state dynamic equations detailed in Section \ref{sec:finitehorizoncontrol} can be solved with numerical toolboxes for finite-horizon optimal control (cf. \cite{wang2009solving} for a MATLAB-based solution). However, for the application of interest, there are major issues which limit the numerical performance, including: 1) model inaccuracies, 2) noisy observations (inaccurate case reports), 3) missing reports (during holidays), 4) unknown or variable parameters (which is inevitable for a highly dynamic complex system, such as a global pandemic), 5) the difficulty of incorporating the start- and end-point boundary conditions from (\ref{eq:transversalityconditions}).

Due to these issues, we propose a novel technique based on optimal state estimation. Accordingly, we integrate the finite-horizon NPI optimizer and the pandemic trend predictor in a classical extended Kalman filter (EKF) and extended Kalman smoother (EKS) \cite{AndersonMoore1979}. Using (\ref{eq:dynamicmodel}), (\ref{eq:observationNewCases}) and (\ref{eq:PontryaginPrinciple}), the state-augmented dynamic equations for the EKF are:
\begin{equation}
\begin{array}{l}
\dot{s}(t) = -\alpha(t)s(t)i(t) + w_s(t)\\
\dot{i}(t) = \alpha(t)s(t)i(t) - \beta i(t) + w_i(t)\\
\dot{\alpha}(t) = -\gamma \alpha(t) + \gamma h[\mathbf{u}^*(t)] + w_{\alpha}(t)\\
\dot{\lambda}_1(t) = [\lambda_1(t) -\lambda_2(t)-(1 - \epsilon)]\alpha(t)i(t) + \eta_{1}(t)\\
\dot{\lambda}_2(t) = [\lambda_1(t) -\lambda_2(t)-(1 - \epsilon)]\alpha(t)s(t) + \beta\lambda_2(t) + \eta_{2}(t)\\
\dot{\lambda}_3(t) = [\lambda_1(t) -\lambda_2(t)-(1 - \epsilon)]s(t)i(t) + \gamma\lambda_3(t) + \eta_{3}(t)\\
n(t)=\alpha(t)s(t)i(t) + v(t)
\end{array}
\label{eq:EKF}
\end{equation}
where the first six equations are the state and co-state dynamics, the last equation is the observation equation and $h[\mathbf{u}^*(t)]$ is the impact of the optimal control calculated from (\ref{eq:optimalinputlinear}) or (\ref{eq:optimalinputquadratic}). The terms $w_s(t)$, $w_i(t)$, $w_\alpha(t)$, $\eta_1(t)$, $\eta_2(t)$ and $\eta_3(t)$ in (\ref{eq:EKF}) represent process noises (due to model inaccuracies), and $v(t)$ is observation noise. Note that for an estimation based on the total number of confirmed cases (instead of the new cases), the last observation equation in (\ref{eq:EKF}) can be replaced with (\ref{eq:observationTotalCases}). 
The dynamic system (\ref{eq:EKF}) may now be numerically solved by using standard EKF and EKS equations. The discretized version of (\ref{eq:EKF}), which is required for the discrete-dime implementation of the EKF and EKS is detailed in the Appendix.

With the above formulation, the finite-horizon optimal NPI prescription and the EKF/EKS-based forecasting problems are unified in a single model. The overall scheme is summarized in Algorithm \ref{alg:prescription}. % and depicted in Fig.\ref{???}.
The MATLAB implementation of this algorithm is available in our online repository \cite{EpidemicModelingCodesSameni2020}. Note however that the prediction and prescription schemes may still be used independently. Specifically, the proposed optimal NPI prescriptor may be combined with any reliable forecasting scheme, including the data-driven LSTM-based models developed during the XPRIZE Pandemic Response Challenge \cite{miikkulainen2021prediction,Lozano2021,Janko2021}.

\begin{algorithm}[tb]
\begin{flushleft}
\caption{Summary of the proposed algorithm
\label{alg:prescription}}
\begin{algorithmic}[1]
\REQUIRE{Historic case reports and NPIs (or an arbitrary scenario from any \textit{predictor model})}
\REQUIRE{The NPI weights $\mathbf{w}(t)$ per region/country}
\REQUIRE The Pareto front tuning parameter $\epsilon\in[0, 1]$. 
\FORALL{Regions}
\STATE Train the compartmental model parameters over historic NPI and case reports (or a predictor model).
\STATE Use EKF and EKS for prediction and prescription of finite-horizon optimal control inputs $\mathbf{u}^*(t)$.
\ENDFOR
\end{algorithmic}
\end{flushleft}
\end{algorithm}

\section{Model training}
\label{sec:modeltraining}
The model parameters $h[\mathbf{u}(t)]$, $\beta$, $\gamma$ and the EKF/EKS parameters require region-wise training or fact-based selection. For this study, we used classical techniques for \textit{Kalman filter engineering}, based on monitoring the properties of the \textit{innovations process} of the Kalman filter to select and automatically adapt the Kalman filter parameters (initial/final states and covariance matrices) over time \cite[Ch. 8]{grewal2007global}. The parameters related to the social and epidemic aspects of the model are explained in the sequel.

\begin{comment}
Add a table containing the initial values of the EKF/EKS parameters.
\end{comment}

The mapping $h[\mathbf{u}(t)]$ was trained over the OxCGRT dataset historic cases. For this, the developed filter was first applied to the historic data by neglecting the explicit relation between the NPI and contact rates (equivalent to $h[\mathbf{u}(t)] = 0$). Referring to (\ref{eq:EKF}), this assumption is equivalent to considering the input-driven fluctuations of $\alpha(t)$ in the process noise $w_\alpha(t)$. Therefore, the entry of the process noise covariance matrix, which corresponds to $w_\alpha(t)$ is increased (as compared with the expected $\alpha(t)$ error) to account for the inaccuracy of the model due to neglecting $h[\mathbf{u}(t)]$. The resulting EKS gives a primary estimate of $\alpha(t)$ over the training period, which in the next stage is given to a constrained LASSO or quadratic polynomial fitter (for the linear and quadratic forms presumed in Section~\ref{sec:NPIhumanInteraction}) to estimate $h[\mathbf{u}(t)]$ using the historic NPI data. Next, the trained model $h[\mathbf{u}(t)]$, together with the historic data is used in a second round of EKS, this time by using the historic NPI and apparently a smaller \textit{a priori} assumption for the variance of $w_\alpha(t)$ (as it no longer accounts for neglecting $h[\mathbf{u}(t)]$). After the secondary EKS, the new estimates of $\alpha(t)$ are once more used to refine the model parameters of $h[\mathbf{u}(t)]$. The refined parameters are stored per country/region for utilization during the prescription phase over real or synthetic scenarios.

The \textit{action-to-effect rate} parameter $\gamma$ was selected intuitively. From various social experiences, it is reasonable to expect a smooth transition in $\alpha(t)$ due to any change in the NPI. This is based on the social experience that imposing any policy on a complex social system is rarely abrupt. Although the transition is region and NPI dependent, in order to reduce the model complexity, we have fixed $\gamma$ to 1/(7\,days)=0.1429\,days$^{-1}$ for all regions/countries.

The recovery parameter $\beta$ was selected by educated guesses from the Centers for Disease Control and Prevention (CDC) reports regarding recovery and contagion periods\footnote{Refer to CDC guidelines for Interim Guidance on Ending Isolation and Precautions for Adults with COVID-19: \url{https://www.cdc.gov/coronavirus/2019-ncov/hcp/duration-isolation.html}}. Accordingly, multiple scientific studies worldwide have reported that an exposed subject is no longer \textit{infectious} after three to four weeks. Of course, this is a stochastic range. To clarify, with an exponential model such as the SI model, in absence of new infected cases ($\alpha = 0$), we find the ratio $i(t_0 + T)/i(t_0)=\exp(-\beta T)$, which can be considered as an exponential law for the \textit{probability of infectiousness after $T$ time-units (days)}. Combining the model with the CDC reports, we derive the following rule of thumb for setting $\beta$:
\begin{equation}
%\beta = \frac{-\log(\text{probability of contagion after $T$ time units})}{T}
\beta = -\log(\text{contagion probability after time $T$})/T
\end{equation}
For the later presented results, we have set the probability of contagion to 0.01 at $T$=21\,days, resulting in $\beta$=0.219\,days$^{-1}$.

Following recent studies \cite{Petersen2020}, the reproduction rate of the pandemic during outbreak was taken to be $\mathcal{R}_0$=2.5 with $\Delta$=1 day, which using (\ref{eq:ReproductionRate}) together with $\beta$ were used to initialize $\alpha(t_0)$, i.e. the contact rate during outbreak. Note that one of the advantages of the EKF/EKS framework is that the model parameters can also be considered as state variables and be \textit{state-augmented} with the other equations to be estimated (or updated over time). This approach can be used for both $\gamma$ and $\beta$ to refine the initial educated guesses.%, as detailed in Section \ref{sec:adaptivecontrol}.  

Finally, the regional/national population sizes, as required for normalizing the total and new contaminated cases to the normalized variables of the SI model were obtained from the United Nations' World Population Prospects 2019 dataset \cite{UNGlobalPopulation}, and was assumed to remain fix over the study, i.e. immigration, inter-border travels, natural birth/deaths have been neglected throughout the study.

\section{Performance monitoring}
\label{sec:performance_monitoring}
Kalman filters have intrinsic mechanisms for monitoring their performance and the consistency of their selected parameters. The so-called \textit{innovation process} is at the heart of performance monitoring. It is defined:
\begin{equation}
	\hat{v}(t) \stackrel{\Delta}{=} x(t) - \hat{x}(t)
	\label{eq:innovation}
\end{equation}
where $\hat{x}(t)$ is an estimate of the desired observation $x(t)$, corresponding to the total number of cases or the new cases, as defined in (\ref{eq:observationTotalCases}) or (\ref{eq:observationNewCases}), respectively. 
For a well-functioning Kalman filter with a single observation (the total number of cases or the new cases), the innovation process has the following properties \cite[Ch 5.3]{AndersonMoore1979}, \cite{GrewalAndrews01}:
\begin{enumerate}[label= P{{\arabic*}})]
\item $\mathbb{E}\{\hat{v}(t)\} = 0$,
\item $\mathbb{E}\{\hat{v}(t) \hat{v}(t') ^ T\} = 0$ for $t \neq t'$,
\item $\gamma(t) \stackrel{\Delta}{=} \mathbb{E}\{\hat{v}(t) \hat{v}(t) ^ T\} = \mathbf{c}(t)^T\bar{\mathbf{P}}(t)\mathbf{c}(t) + r(t)$,
\end{enumerate}
where $\bar{\mathbf{P}}(t)$ is the covariance matrix of the state vector estimation error before the $t$\textsuperscript{th} measurement, $\mathbf{c}(t)$ is the state-to-observation map Jacobian, and $r(t)$ is the observation noise variance (defined in Appendix \ref{sec:EKFEKS}). P1 and P2 guarantee that the process noise is zero-mean and spectrally white, and P3 assures that the Kalman filter's estimate of the observation noise is in accordance with the presumed model parameters. The violation of any of these properties is an indication of parameter mis-selection, e.g., the state or observation noise covariance/variance. Although the above properties originally belong to the linear Kalman filter, they can be equally used to monitor the well-functioning of the EKF/EKS. In other contexts, indexes have been proposed based on the above properties for tracking the performance of Kalman filters and the selection of their parameters \cite{GrewalAndrews01}, \cite[Ch. 6]{AndersonMoore1979},  \cite{JamshidianSameniJutten2019}. These monitoring indexes are provided in the open-source implementation of the hereby developed EKF \cite{EpidemicModelingCodesSameni2020}, and are useful for tuning the EKF/EKS model parameters per region/country.

\section{Results}
\label{sec:results}
\subsection{Forecasting}
\begin{figure}[tb]
\includegraphics[trim=0.9in 0in 0.9in 0.0in,clip,width=0.9\columnwidth]{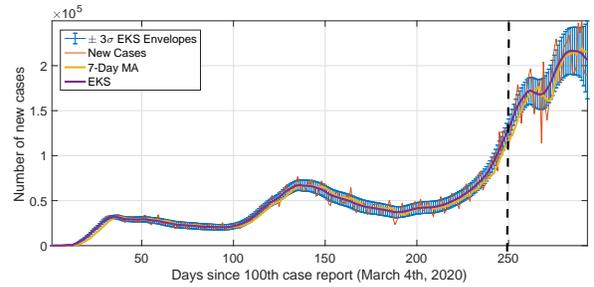}
\caption{New cases trend tracking using the proposed EKS on daily reported cases of the US, since Mar 4\textsuperscript{th} 2020 (the 100\textsuperscript{th} case report). The raw noisy daily reports have been adopted from the OxCGRT dataset \cite{OxCGRT2020}. The smoothing period is up to day 250 and used to forecast the trends thereafter.}
\label{fig:EKSSampleResult}	
\end{figure}

Fig.~\ref{fig:EKSSampleResult} shows the smoothing result of the number of new cases using the proposed EKS on daily reported cases of the US from March 4\textsuperscript{th} 2020 (the 100\textsuperscript{th} case report date), until November 9\textsuperscript{th} 2020. These 250 days of smoothing is followed by 40 days of forecasting by using an EKF. The $\pm 3 \sigma$ standard deviation envelopes obtained from the error covariance matrix estimates of the EKF are also plotted on the trend estimates. It is seen that as we estimate the number of cases farther in the future, the confidence intervals enlarge, i.e. the forecasts become less reliable.

\begin{figure}[tb]
\includegraphics[trim=0.3in 0in 0in 0in,clip,width=\columnwidth]{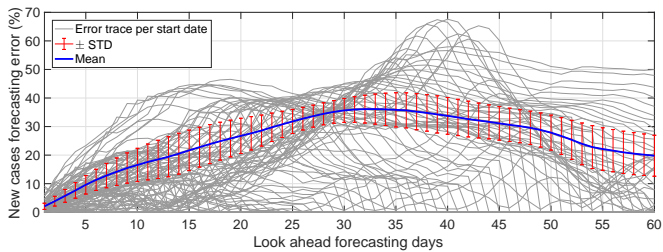}
\caption{New cases look-ahead forecasting error percentage for the US, with variable start dates from from Dec 3\textsuperscript{rd} 2020 to Mar 4\textsuperscript{th} 2021. Real data and NPI from Dec 3\textsuperscript{rd} 2020 to Mar 4\textsuperscript{th} 2021 have been used for model training.}
\label{fig:ForecastingErrors}	
\end{figure}

In order to show the forecasting accuracy, another experiment was conducted, in which the EKF/EKS NPI prescription model was trained by the real US new cases and NPI data from March 4\textsuperscript{th} 2020 through December 3\textsuperscript{rd} 2020. The trained model was applied to the ground truth data from December 4\textsuperscript{th} 2020 to March 4\textsuperscript{th} 2021, to forecast the number of new cases from 1 to 60 days ahead, having only the actual NPIs as the input. Fig.~\ref{fig:ForecastingErrors} shows the percentage of forecasting error for different start dates from December 3\textsuperscript{rd} 2020 to March 4\textsuperscript{th} 2021, versus the number of look-ahead forecasting days. As expected, the model is generally more accurate for short-term forecasting as compared to long-term forecasting. Note that although we can see that the look-ahead forecasting errors have slightly decreased beyond 30 days, the observation is not generalizable and may not be associated to the accuracy of the forecasting model, as it highly depends on the dynamics of the pandemic, ground-truth data, NPI policies adopted by a country, and many other socioeconomic factors. We further elaborate on this point in the discussion.

\subsection{NPI prescription}
\begin{figure*}[tb]
    \centering
\subfigure[Afghanistan]{\label{fig:Afghanistan}\includegraphics[trim=0.8in 0in 1.1in 0in,clip,width=.23\textwidth]{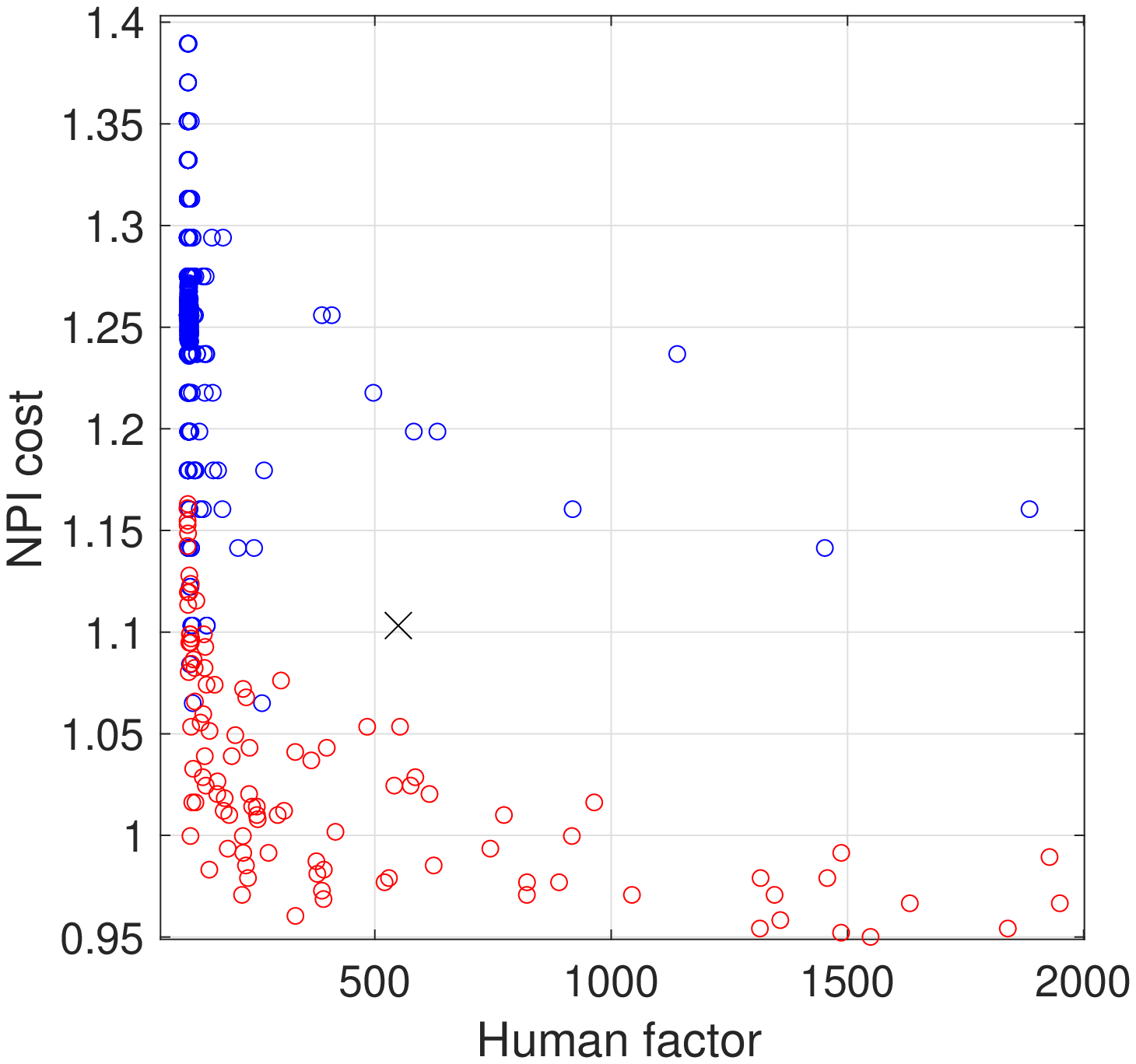}}
\subfigure[Argentina]{\label{fig:Argentina}\includegraphics[trim=0.8in 0in 1.1in 0in,clip, width=.23\textwidth]{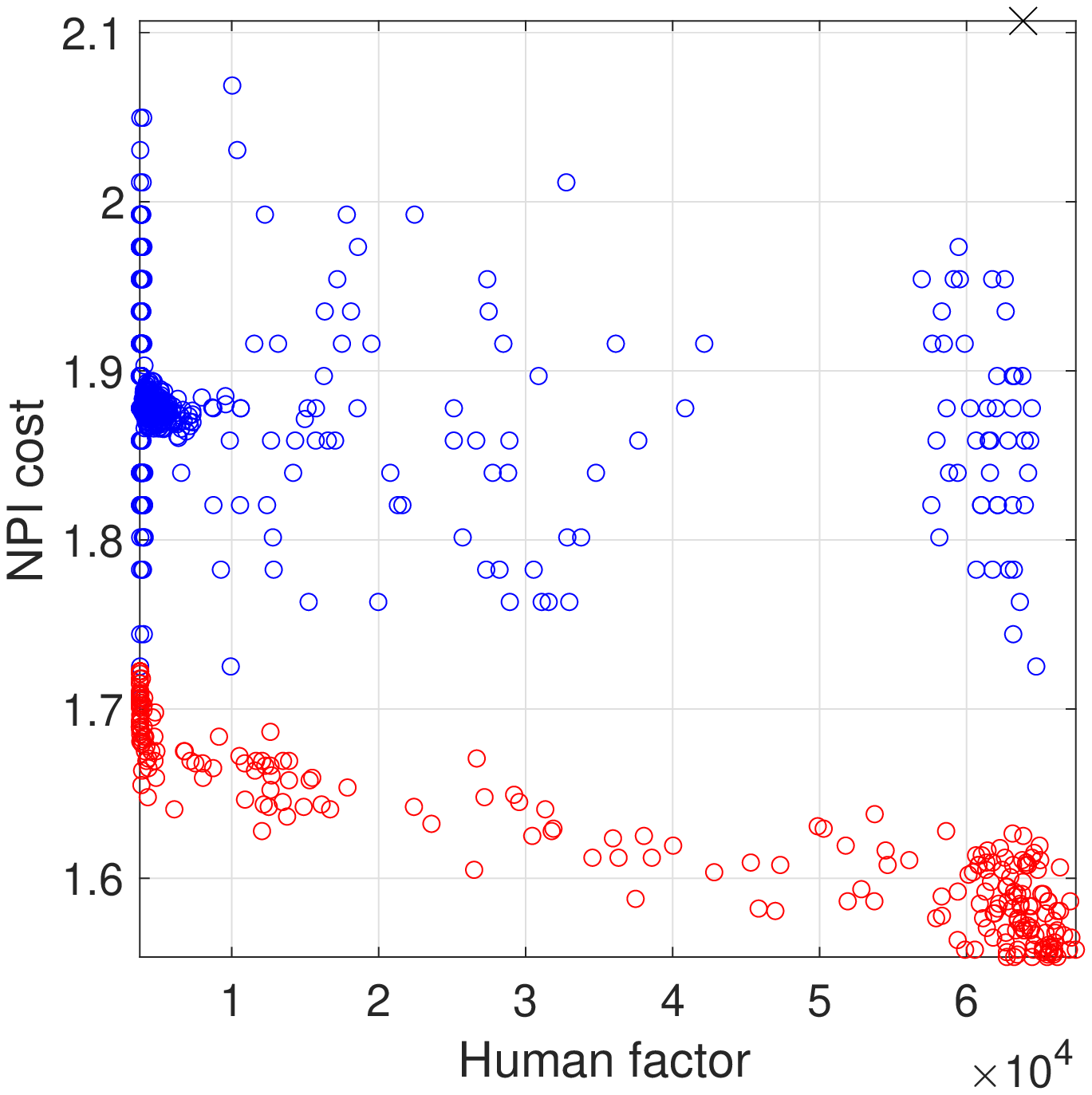}}
\subfigure[Brazil]{\label{fig:Brazil}\includegraphics[trim=0.8in 0in 1.1in 0in,clip, width=.23\textwidth]{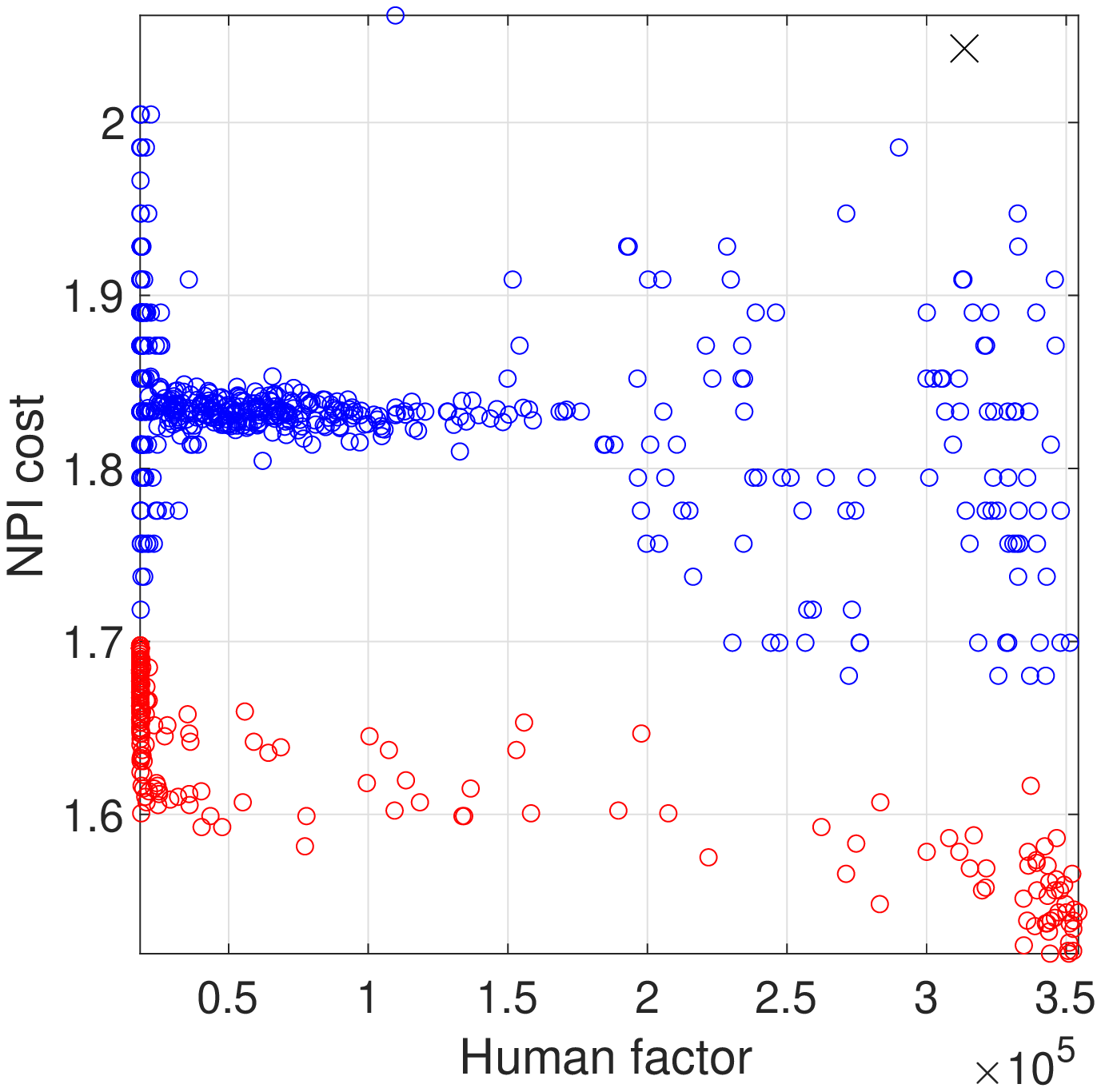}}
\subfigure[China]{\label{fig:China}\includegraphics[trim=0.8in 0in 1.1in 0in,clip, width=.23\textwidth]{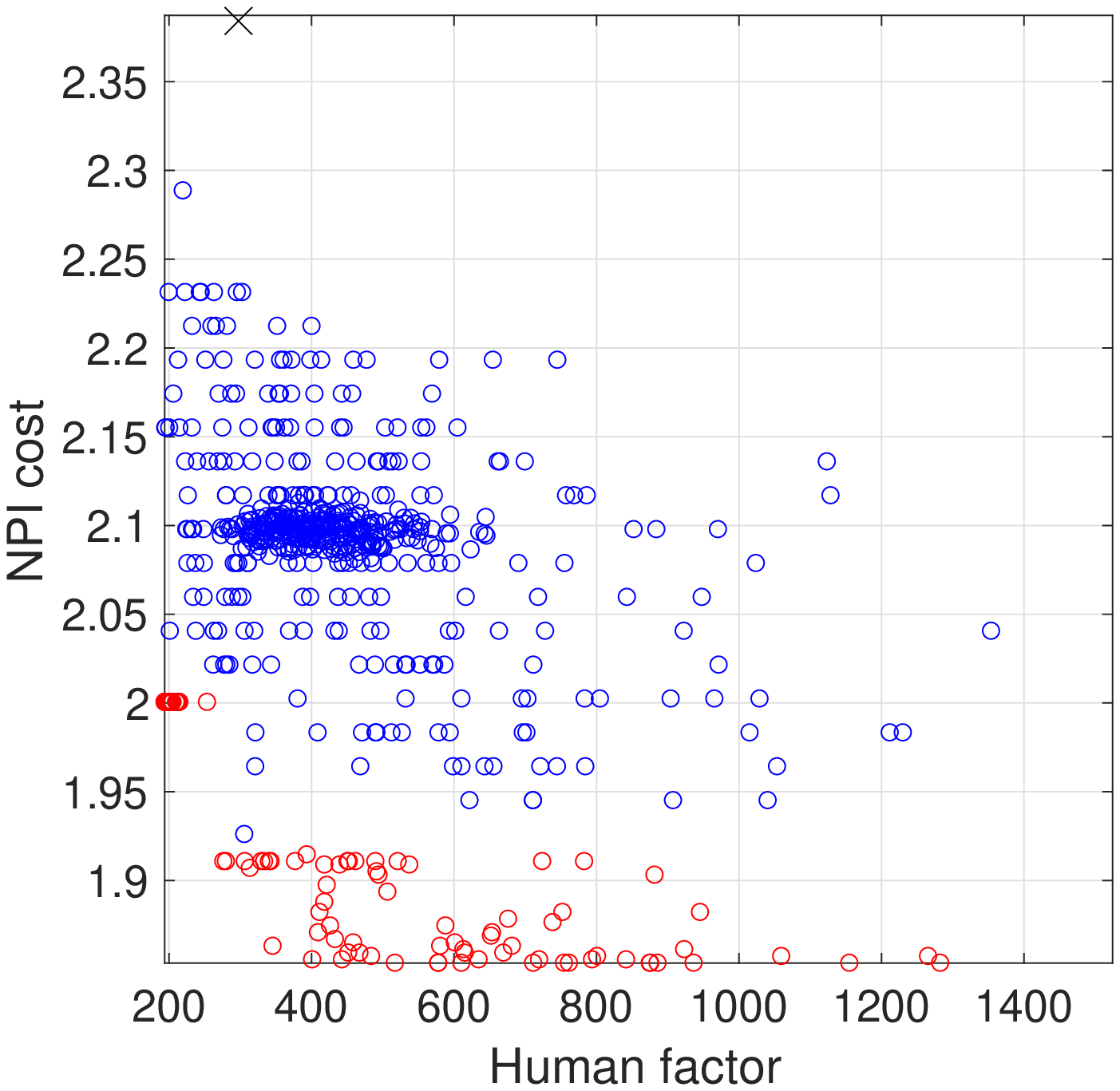}}
\subfigure[France]{\label{fig:France}\includegraphics[trim=0.8in 0in 1.1in 0in,clip, width=.23\textwidth]{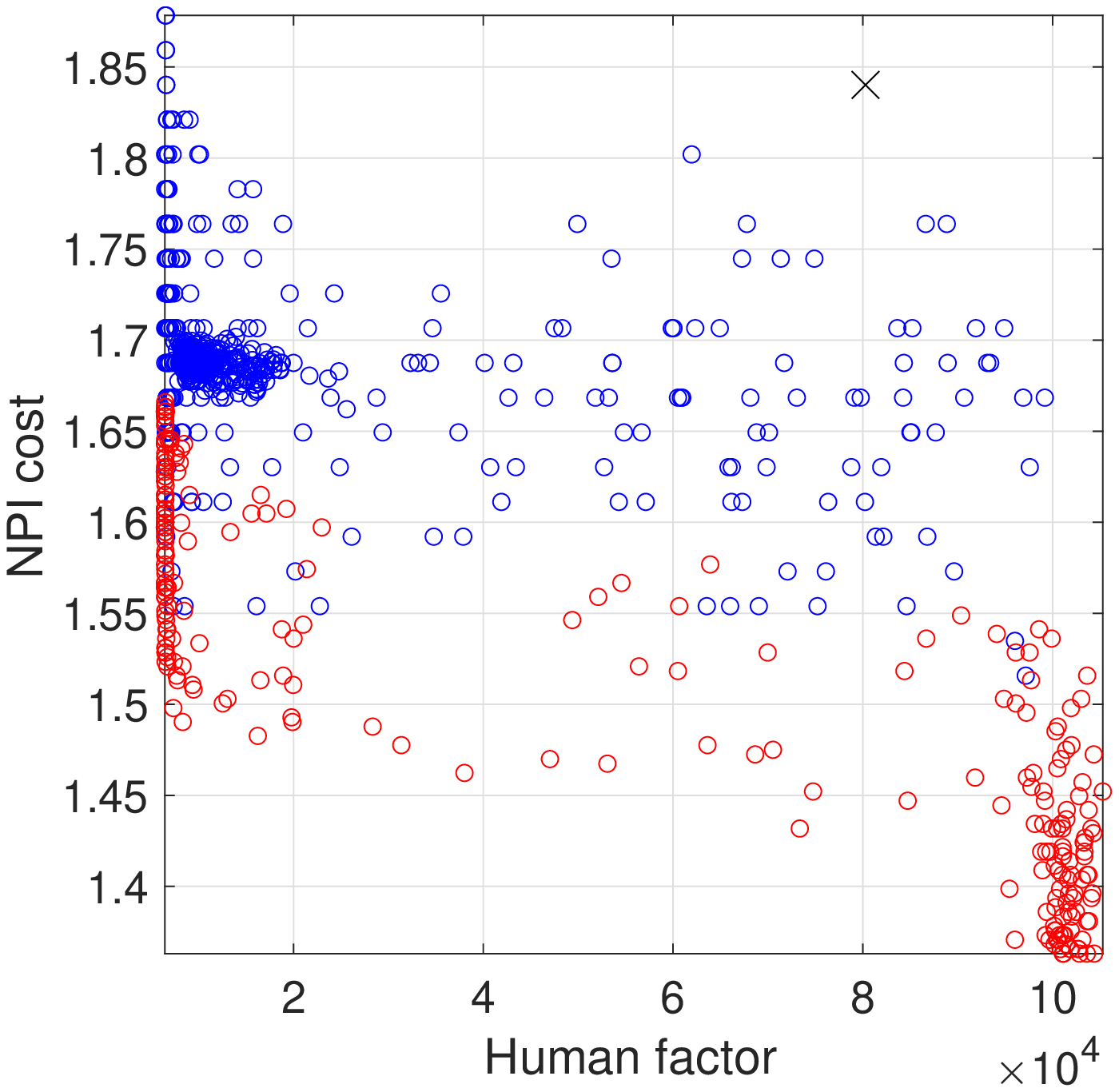}}
\subfigure[Iran]{\label{fig:Iran}\includegraphics[trim=0.8in 0in 1.1in 0in,clip, width=.23\textwidth]{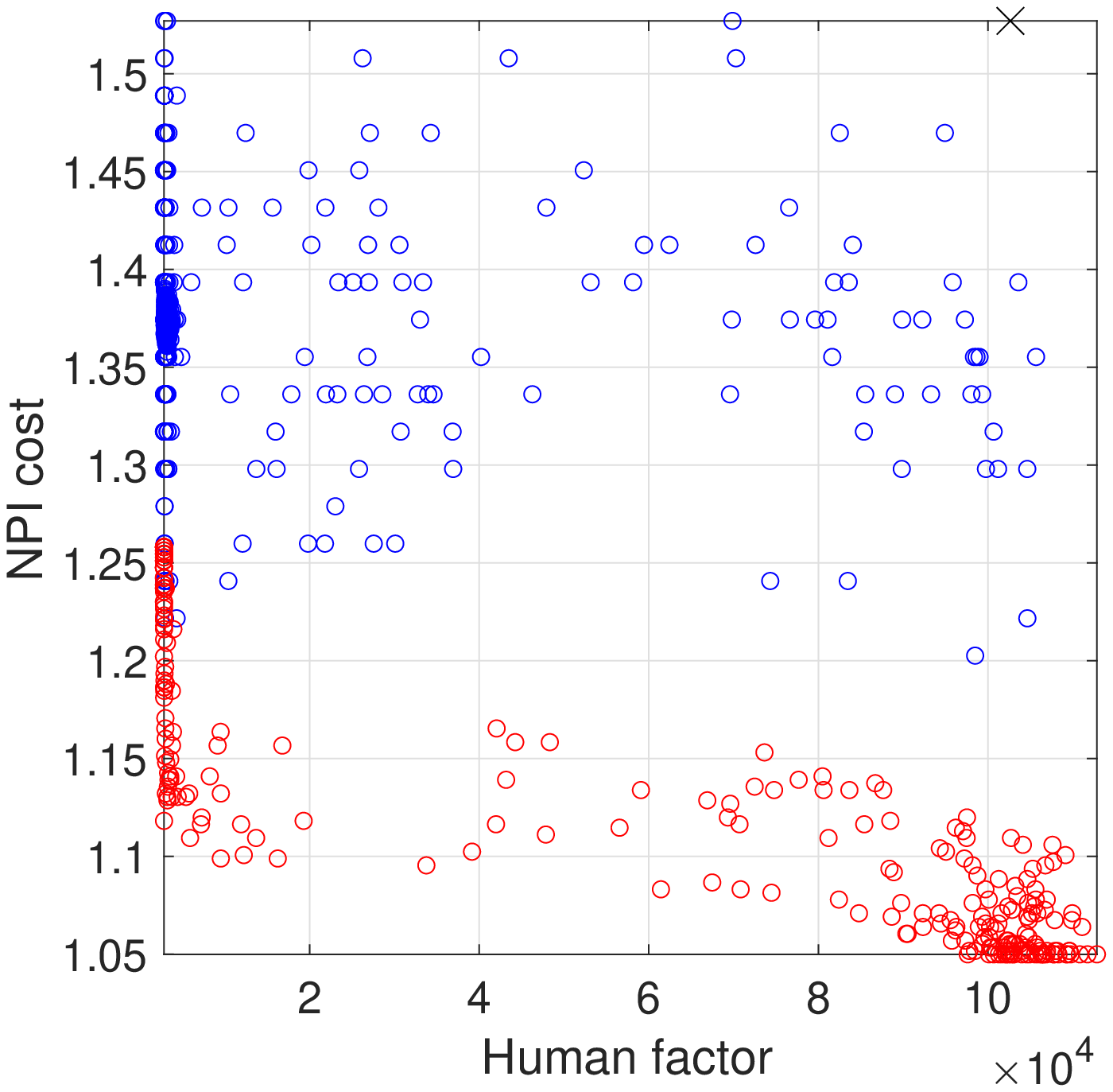}}
\subfigure[Italy]{\label{fig:Italy}\includegraphics[trim=0.8in 0in 1.1in 0in,clip, width=.23\textwidth]{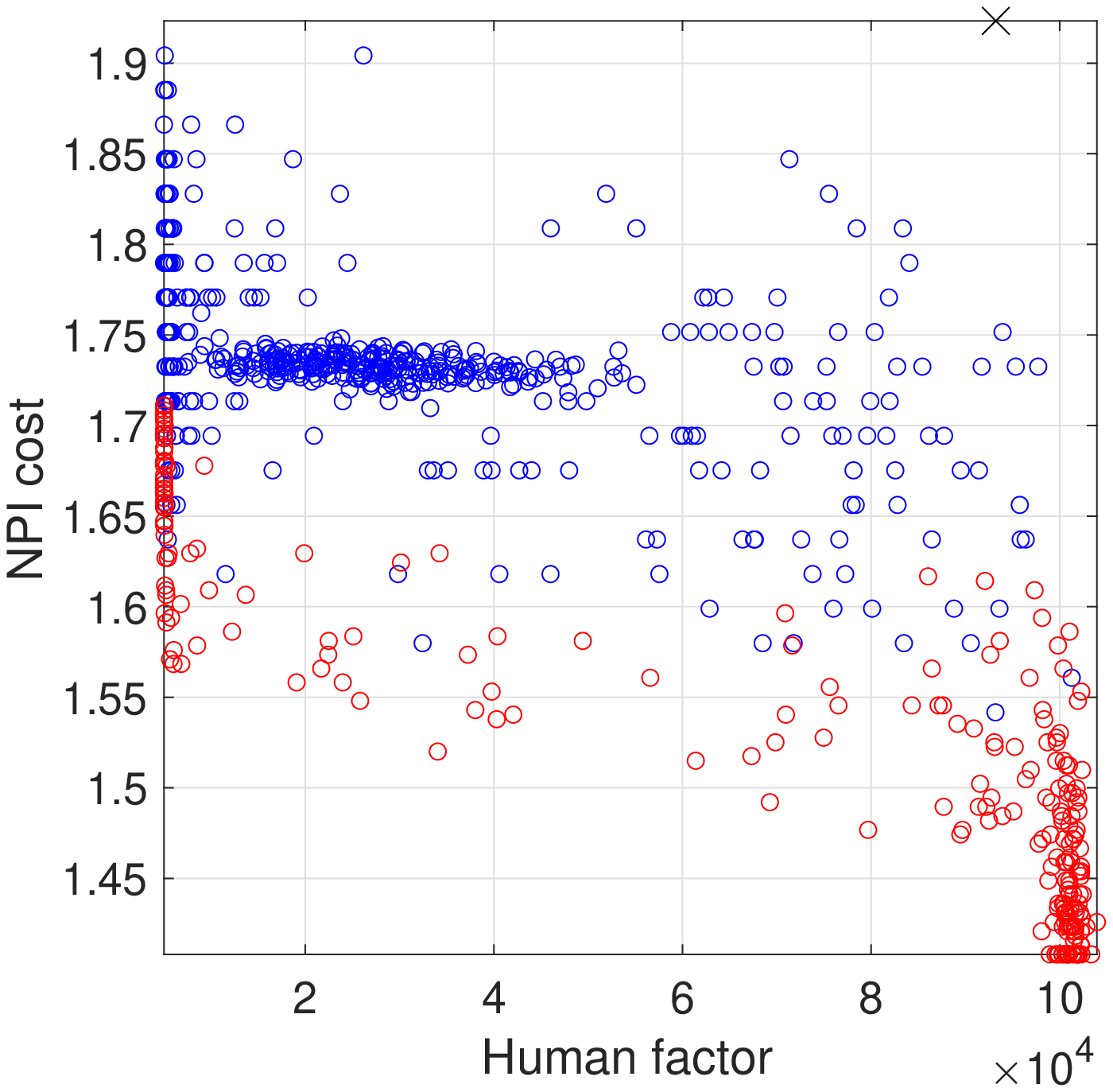}}
\subfigure[Germany]{\label{fig:Germany}\includegraphics[trim=0.8in 0in 1.1in 0in,clip, width=.23\textwidth]{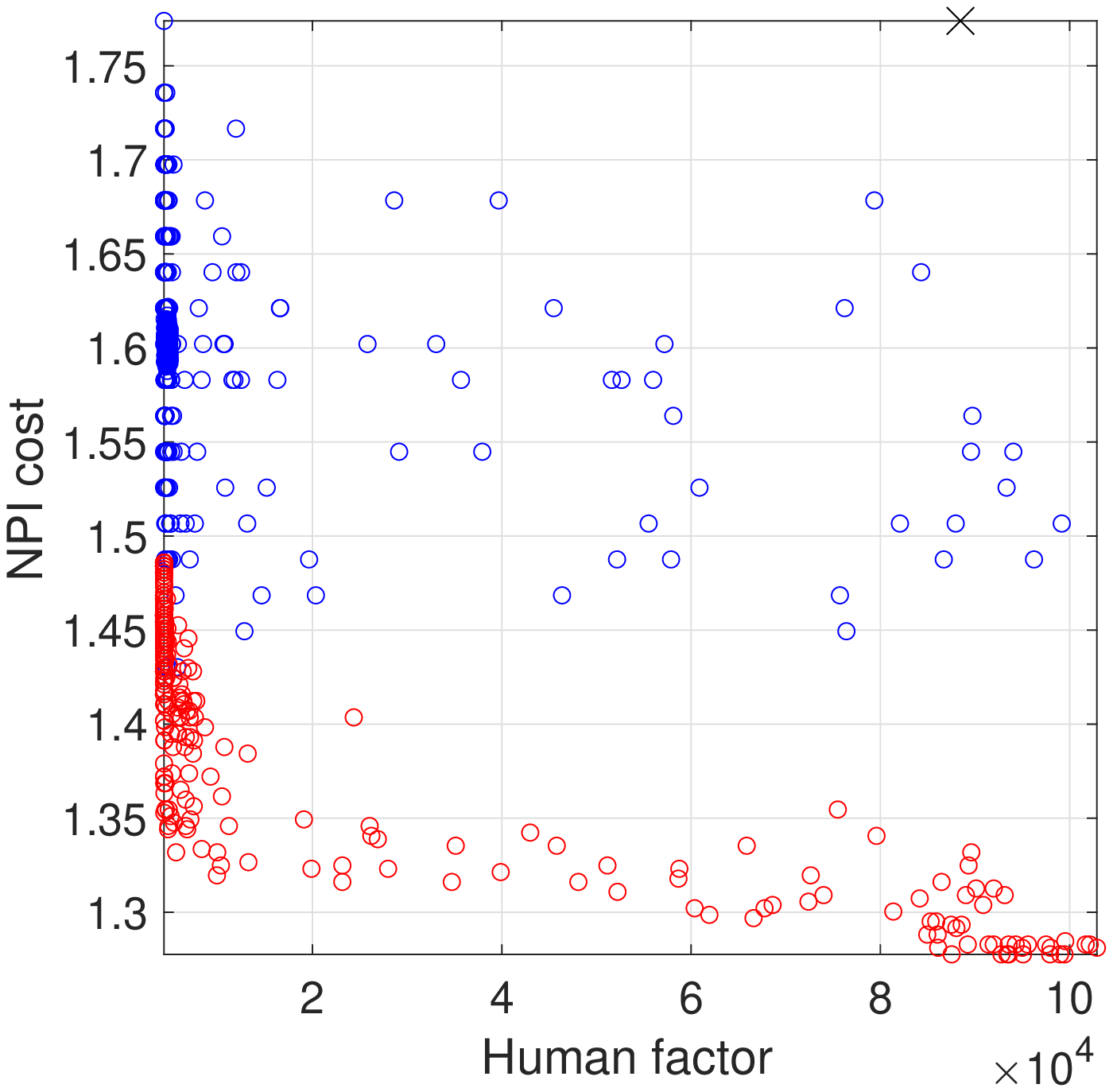}}
\subfigure[South Africa]{\label{fig:SouthAfrica}\includegraphics[trim=0.8in 0in 1.1in 0in,clip, width=.23\textwidth]{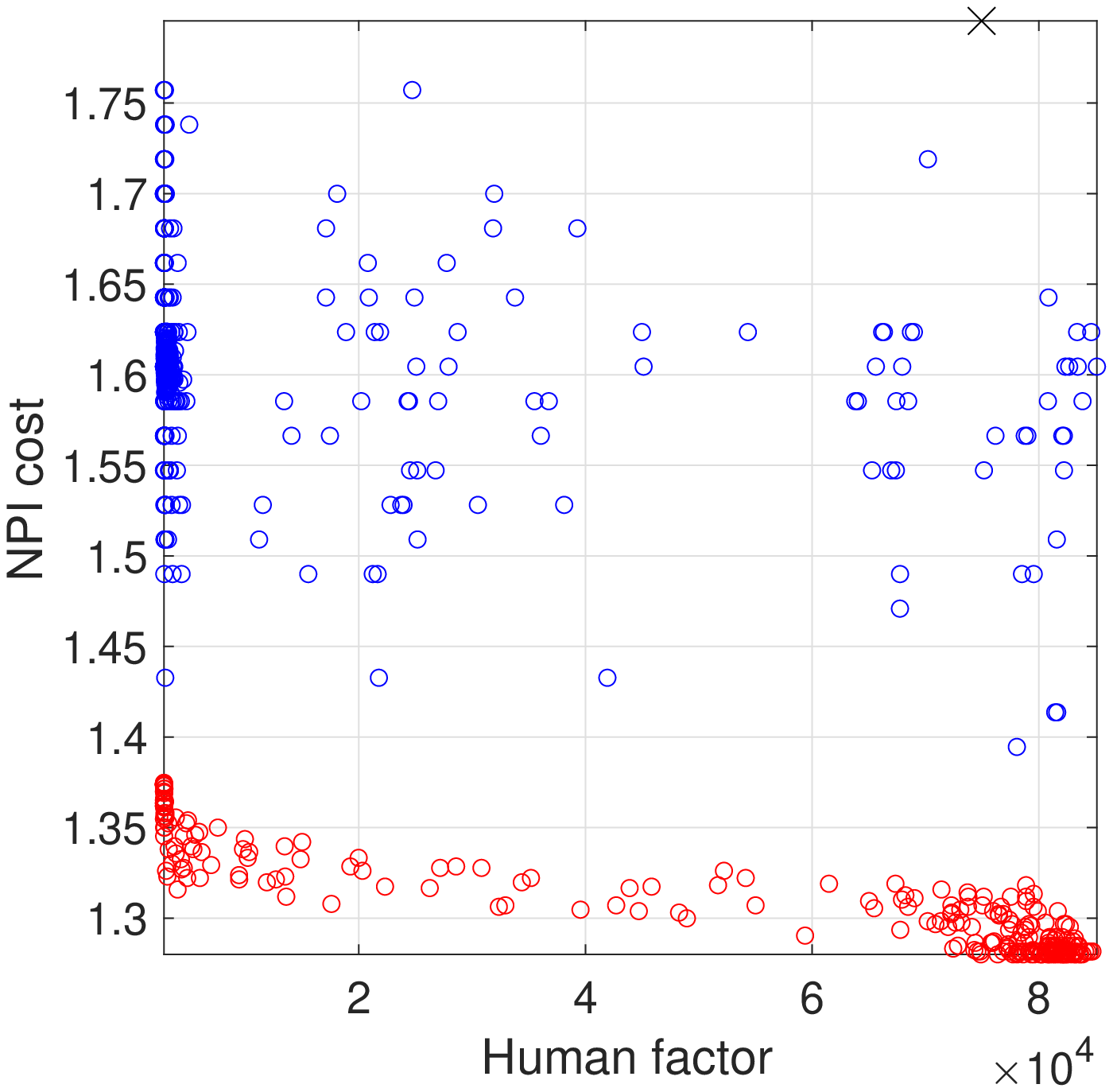}}
\subfigure[Spain]{\label{fig:Spain}\includegraphics[trim=0.8in 0in 1.1in 0in,clip, width=.23\textwidth]{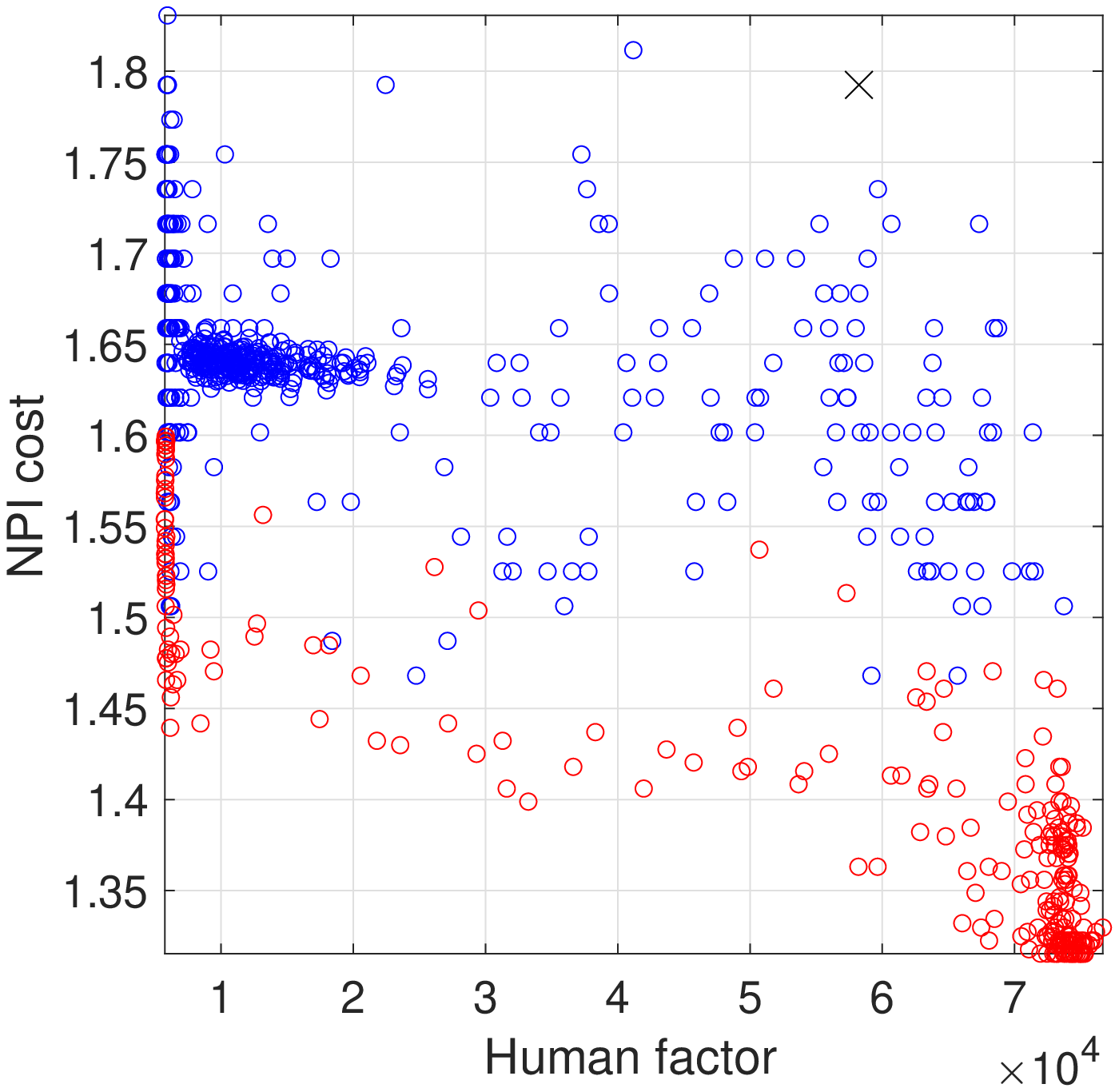}}
\subfigure[UK]{\label{fig:UK}\includegraphics[trim=0.8in 0in 1.1in 0in,clip, width=.23\textwidth]{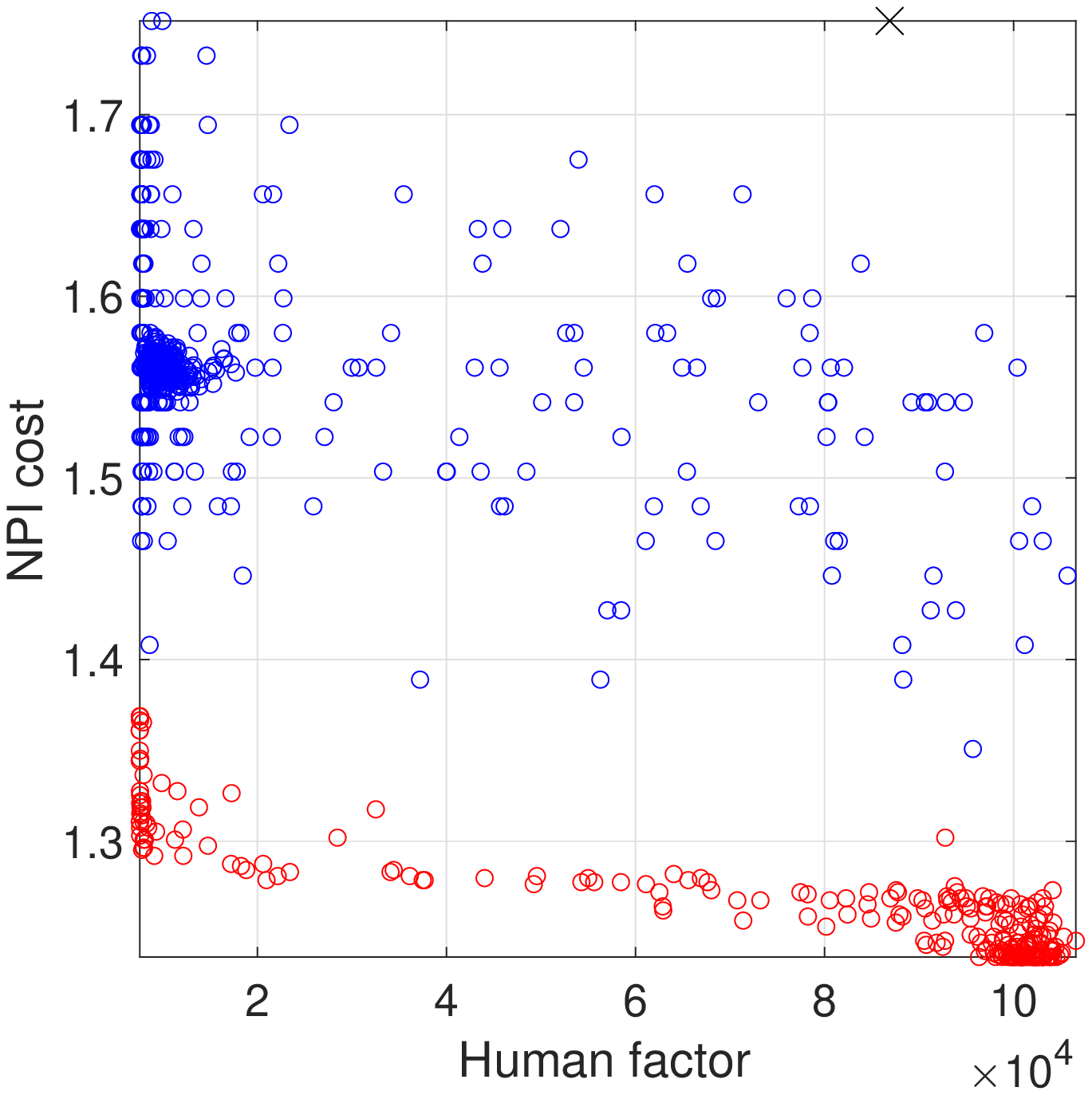}}
\subfigure[US]{\label{fig:US}\includegraphics[trim=0.8in 0in 1.1in 0in,clip, width=.23\textwidth]{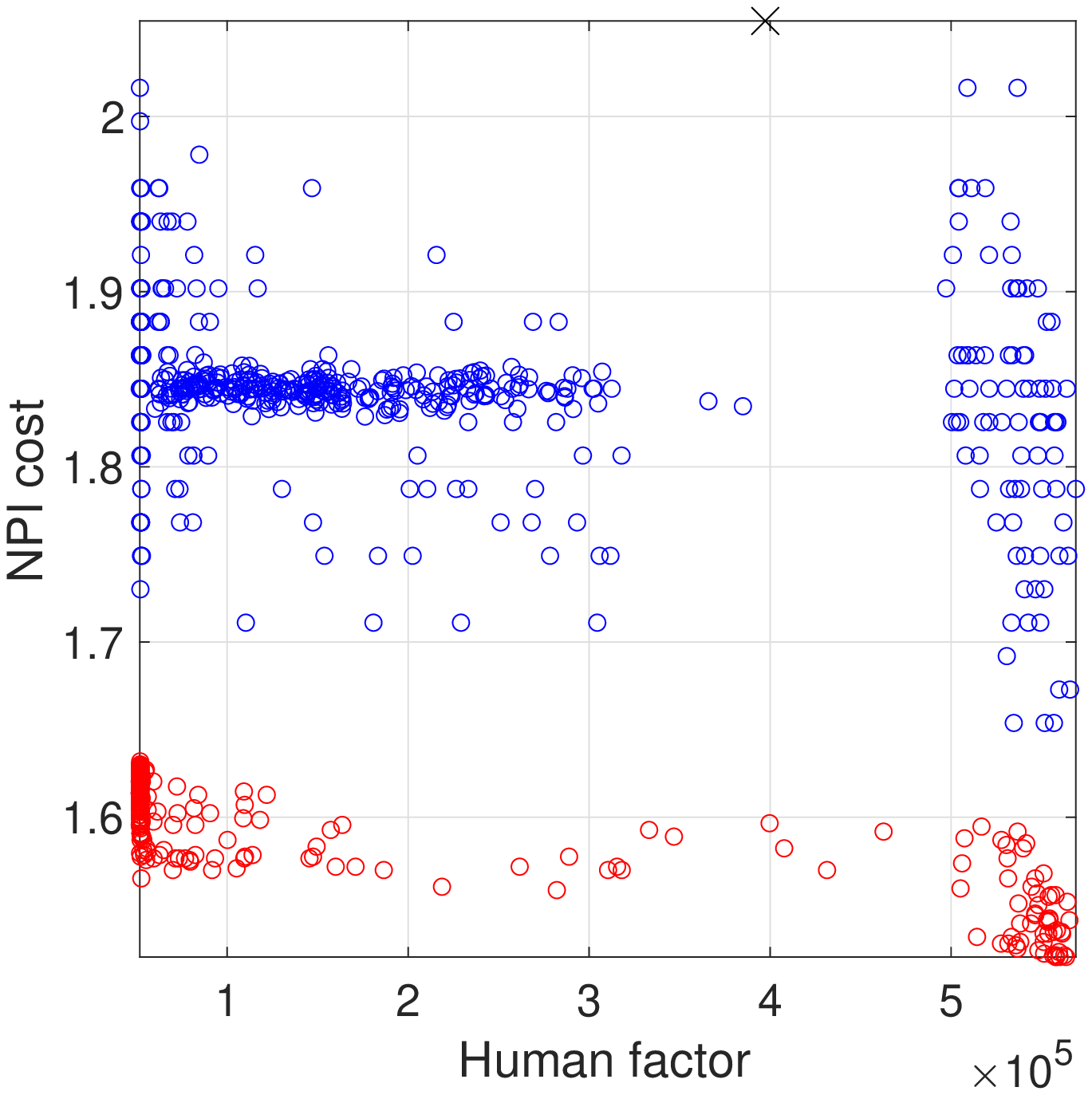}\label{fig:USnpiBiObjective}}
\caption{Biobjective optimization space for sample countries. Black cross: fixed NPI (continuing current policies); Blue: random NPI inputs (both constant and variable over time); Red: finite-horizon optimal input for 250 $\epsilon\in[0, 1]$. $h[\mathbf{u}(t)]$ was found by linear regression over historic NPI from Jan 1\textsuperscript{st} 2020 to Feb 7\textsuperscript{th} 2021, using a LASSO with positive coefficients constraint.}
\label{fig:paretofronts}
\end{figure*}

The OxCGRT dataset has above 300 countries and regions (states). Due to inconsistencies in the reported COVID-19 cases, some of the countries/regions were omitted from the study during the XPRIZE Challenge and the proposed prediction-prescription algorithm was trained and applied to a total number of 235 regions/countries, with arbitrary NPI cost weights. The training period for the model was from January 1\textsuperscript{st} 2020 to Feb 7\textsuperscript{th} 2021, and the test phase was from February 8\textsuperscript{th} 2021 up to May 7\textsuperscript{th} 2021 (the preparation date of the manuscript). As proof of concept, the bi-objective optimization space of the NPI cost $J_1$ vs the human factor cost $J_0$ are shown in Fig.~\ref{fig:paretofronts}, for several countries worldwide. Accordingly, each point in this figure corresponds to a $(J_0, J_1)$ pair for a sequence of NPI scenarios over the test phase. In Fig.~\ref{fig:paretofronts}, the red points are the result of the proposed method for different values of $\epsilon$ (the bi-objective optimization free parameter). The black crosses correspond to the scenario of continuing the latest NPI policy of each government at the end of the training date, up to the end of the testing date. Finally the blue points correspond to a pool of \textit{random constant} stringencies $\mathbf{u}(t) = \bm{\kappa}$, $\bm{\kappa}\sim U[\mathbf{u}^{\min}, \mathbf{u}^{\max}]$, and \textit{random variable} stringencies $\mathbf{u}(t)\sim U[\mathbf{u}^{\min}, \mathbf{u}^{\max}]$, i.e. uniformly distributed between $\mathbf{u}^{\min}$ and $\mathbf{u}^{\max}$. In all cases, the user defined NPI weight vector was chosen to be equal for all NPI ($\mathbf{w}(t) = \mathbf{1}$), i.e., the NPI were considered equally important for the policymaker.

As a bi-objective problem, the Pareto efficient front comprises of the NPI points which either have a smaller value of $J_0$ or $J_1$, while the non-efficient solutions are the ones for which there exists at least a point that gives a smaller cost of both $J_0$ and $J_1$. In other words, a Pareto efficient solution should outperform any other solution either in its cost or efficiency. As trivial cases, the maximum stringency case $\epsilon = 0$ (maximal enforcement of social limitations, to minimize human cost) and the minimum stringency case $\epsilon = 1$ (no social constraints, to minimize costs of intervention) are both Pareto efficient; the former minimizes the human losses and the latter minimizes the socioeconomic cost of intervention. Apparently, policymakers prefer a balance between these two objectives. Therefore, the Pareto efficient NPI policies should be close to the origin or along one of the left or bottom axes. From Fig.~\ref{fig:paretofronts}, it is seen that none of the NPI policies adopted by countries/regions were optimal (assuming equal NPI weights), and despite the significant disparities between the studied countries, the Pareto optimal points all belong to the proposed algorithm. Note that the ``optimal point'' is clearly a function of the bi-objective parameter $\epsilon$, eventually selected by policymakers.

\begin{figure}[tb]
\includegraphics[trim=0.8in 0.05in 1.1in 0.3in,clip,width=0.93\columnwidth]{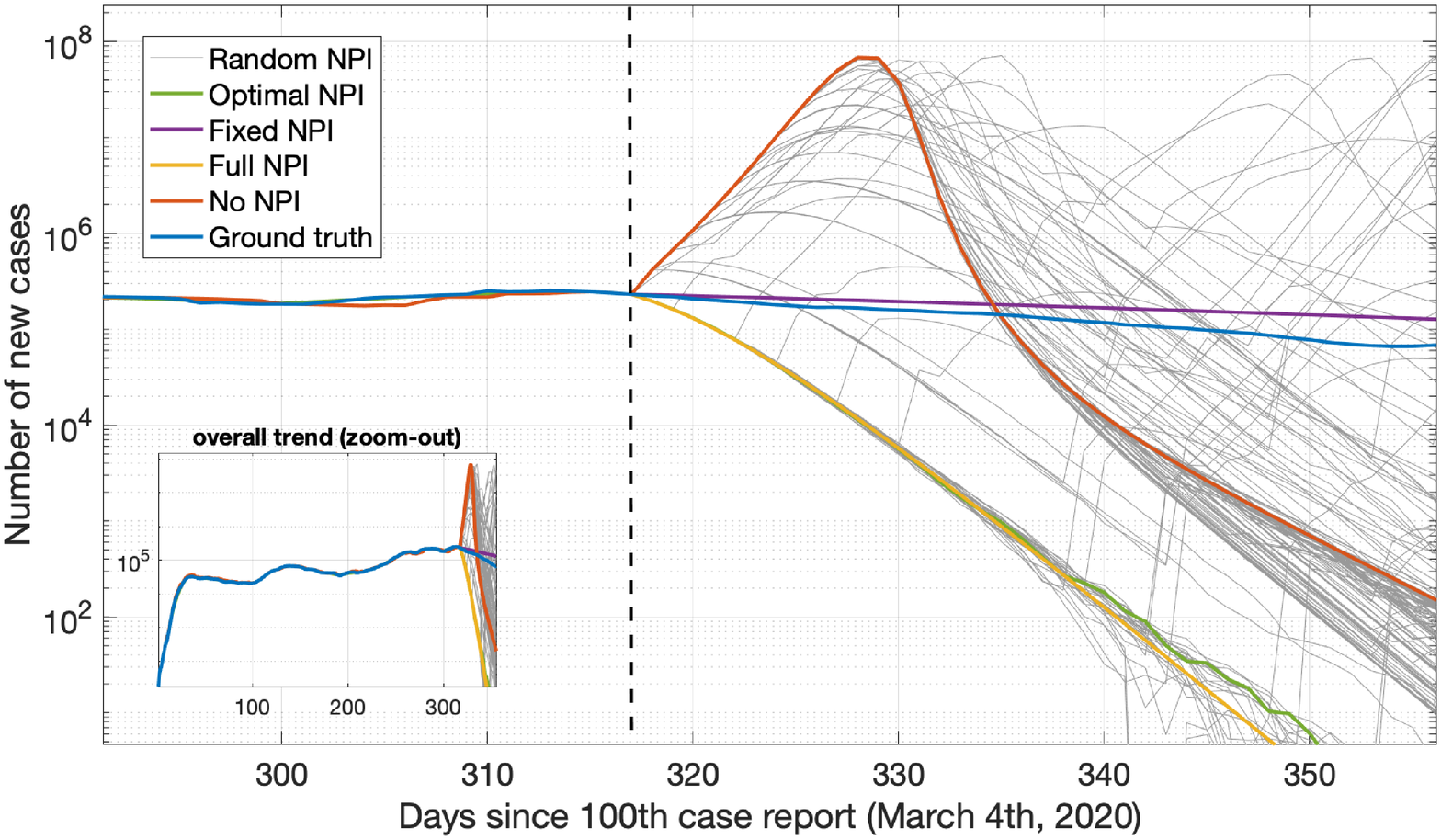}
\caption{The number of new case estimates in the US for different NPI scenarios: optimal NPIs corresponding to 500 random $\epsilon \in[0, 1]$, fixed NPI (no changes in the last training date NPI), full NPI (maximum stringency), no NPI (no stringency), and optimal NPI (numerically found to be $\epsilon = 8.235\times10^{-8}$ corresponding to the NPI scenario with closest normalized distance from the bi-objective origin in Fig.~\ref{fig:USnpiBiObjective}, which is a compromise between NPI cost and human factors), compared to the ground-truth (official reports). Seven-day averages of official reports from Mar 4\textsuperscript{th} 2020 (the US 100\textsuperscript{th} case report date) through Jan 15\textsuperscript{th} 2021 have been used for NPI model training and applied to ground-truth data until Mar 15\textsuperscript{th} 2021. Notice the scenarios that result in the catastrophic extreme case of total population contamination (\textit{herd immunity} by infection) in a short period after NPI removal.}
\label{fig:US_NPI_Scenarios}	
\end{figure}

A similar performance was obtained for all the 235 studied countries and regions. As a case study, Fig.~{\ref{fig:US_NPI_Scenarios}} demonstrates the estimates number of new case in the US for different NPI scenarios, by using seven-day averages of official reports from March 4\textsuperscript{th} 2020 (the US 100\textsuperscript{th} case report date) through January 15\textsuperscript{th} 2021 for NPI model training. The trained model was applied to the ground-truth data from January 16\textsuperscript{th} 2021 to March 15\textsuperscript{th} 2021 for test. Apparently, the hypothesized NPI scenarios may only be evaluated by simulation using our trained EKF forecasting model or other models. The compared scenarios are: A) optimal NPIs corresponding to 500 random $\epsilon \in[0, 1]$, B) fixed NPI (no changes in the NPI after January 15\textsuperscript{th} 2021), C) maximum stringency, D) zero stringency, and E) a compromised NPI on the bi-objective space Pareto front corresponding to the NPI scenario with the smallest normalized distance from the US bi-objective plane origin in Fig.~\ref{fig:USnpiBiObjective}. The latter point is a compromise between NPI cost and human factors, and was numerically found to correspond to $\epsilon = 8.235\times10^{-8}$. In Fig.~\ref{fig:US_NPI_Scenarios}, the ground-truth new number of cases (the result of the actual US NPI policy over the test period) is also shown for comparison. Accordingly, the human cost of the last scenario, which is Pareto efficient is very close to the full-stringency case, without imposing the maximum stringency socioeconomic cost.

\subsection{Processing load}
Since the Pareto front solutions are found by mathematical derivation (rather than trial and error or cumbersome searches), the proposed framework is extremely computationally efficient and the run-time for testing arbitrary scenarios is minimal. The MATLAB version of the codes applied to all regions and countries (235 in total), takes less than 30\,s to train over the historic cases, on a MacBook~Pro laptop with 2.3\,GHz Quad-Core Intel Core\,i7 and 32\,GB of memory, without notable optimizations. The run-time on the test scenarios takes about 15\,s in total for all regions, as it contains only one EKS stage during the test period, per region/country.

\section{Discussion}
\label{sec:discussion}
The proposed algorithm for predicting pandemic trends and prescribing Pareto optimal NPI policies has multiple advantages over mere data-driven machine learning algorithms. The highlights of this algorithm include: 
    \begin{itemize}[leftmargin=*]
    \item The method is based on theoretical derivations and within the scope of the proposed compartmental model accuracy (which is asymptotically accurate for region/country-level population sizes), it gives accurate \textit{Pareto efficient} solutions.
    \item The operation point on the Pareto front can be tuned by a single parameter $\epsilon\in[0, 1]$ selected by the policymaker, where the corner case $\epsilon=0$ neglects the NPI cost (in favor of the human factor) and $\epsilon=1$ neglects the human factor (in favor of intervention cost).
    \item The prediction and prescription problems are integrated in a unified framework. Nevertheless, the method is applicable to both real-world data and any other machine-learning based technique, which accurately predicts pandemic trends from historic data (see for example \cite{miikkulainen2021prediction}).
    \item This framework can be used for \textit{targeted pandemic control}, where strategists can target specific infection bounds that match the medical resources of a country/region, over a fixed or maximally bounded period of time. Therefore, apart from the optimal NPI and fatality rate objectives, such scenarios can also be considered: \textit{``how to bring the pandemic reproduction rate below 0.8 by a specific time?''}, or \textit{``how to bring the new cases below 200 per day in less than two months?''} The training phase of the pandemic over historic data together with the forecasting model can be used to study the feasibility of such scenarios and the prescription of the required NPI policies that would achieve such goals.
    \item Both the model parameters and NPI cost weights can be updated over time. Therefore, unprecedented events such as vaccination or virus mutation effects can be integrated in the model with appropriate training. According to the so-called \textit{principle of optimality}, ``any portion of an optimal control trajectory is optimal \cite[Sec 6.4]{naidu2003optimal},'' which implies that optimality of future actions is independent of the past. Therefore, the prescribed optimal control strategy may be adopted at any point, regardless of the past actions of a region/nation. 
    \item The computational efficiency of this algorithm permits its combination with other machine learning methods to reduce the search space and to improve the accuracy on other datasets and under more complicated models such as the \textit{long short-term memory} (LSTM), as in \cite{miikkulainen2021prediction}. This feature is specifically useful for training the NPI to inter-human contact map $h(\cdot)$, which requires learning.
    \item The predictor part of the model gives confidence intervals during both the prediction and prescription steps of the algorithm. Therefore, the performance and well-function of the algorithm can be continuously monitored and adapted by using the innovation process, as detailed in Section~\ref{sec:performance_monitoring}.
    \item The proposed framework is extendable to pharmaceutical intervention plans and vaccinations, whenever sufficient data is available to design and train alternative compartmental models for these factors.
\end{itemize}

As a point of reservation, we should note that the scope of algorithmic and machine learning-based pandemic trend forecasts and prescriptions should not be overestimated or exaggerated. As an extreme case, consider a forecasting model which predicts that the daily new cases will reach below a certain threshold in several years. We can debate that such long-term speculations are neither scientific nor of any practical use. Because the intrinsic dynamics of pandemics (with or without NPIs, vaccination or herd immunity), guarantee that at some point in the future, there will no longer be any susceptible groups (neglecting the virus mutations). It is therefore essential to demonstrate that any forecasting is fact-based, nontrivial, and predictable from previous observations (in the causal sense).%%Similarly, in global data challenges where hundreds or thousands of groups compete}

\section{Conclusion and future work}
In this research, a model-based approach was used for the prediction and prescription of NPIs that best balance between an arbitrary weighted-cost of interventions and the human factors (number of new cases) during a pandemic. The proposed algorithm and the prescribed NPIs were proved to be Pareto optimal, within the scope of the utilized compartmental model accuracy. Software implementations of the proposed algorithms are online available at \cite{EpidemicModelingCodesSameni2020}. In future studies, different aspects of this framework can be extended and improved. Specifically, advanced ML algorithms can be used for training the NPI to contact rate function $\mathbf{h}(\cdot)$. The LSTM is specifically a promising approach for this purpose.

For regions which have access to additional data (e.g., the number of hospitalized, number of vaccinated, fatality rate of the virus, the population age pyramid, etc.), more accurate models such as the fatal \textit{susceptible-exposed-infected-recovered} (SEIR) can be used \cite{sameni2020mathematical}. Other indexes such as the daily death reports can be augmented as additional observation equations in the dynamic model (\ref{eq:EKF}), and will help to increase the EKF/EKS accuracies. Theoretical aspects of the proposed EKF/EKS frameworks, including stability conditions, parameter identifiability and robustness to parameter and modeling errors also require further studies.

\appendices

\section{The discrete-time model\label{sec:discretized}}
%\appendix[The discrete-time model\label{sec:discretized}]
For a discrete-time implementation of the EKF/EKS, the discrete form of the dynamic system (\ref{eq:EKF}) is required. Accordingly, we define the discrete variables
\begin{equation*}
    \begin{array}{l}
         \!\!\bm{s}_k\!=\![s(k\Delta), i(k\Delta), \alpha(k\Delta)]^T  \\
         \!\!\bm{w}_k\!=\![w_s(k\Delta), w_i(k\Delta), w_\alpha(k\Delta), \eta_1(k\Delta), \eta_2(k\Delta), \eta_3(k\Delta)]^T\\
         \!\!n_k\!=\!n(k\Delta), \quad c_k\!=\!c(k\Delta), \quad v_k\!=\!v(k\Delta)\\
    \end{array}
\end{equation*}
where $\Delta$ is the discretization time unit. Assuming that $\Delta$ is small as compared with the variations of the pandemic trends, a first order discrete approximation of (\ref{eq:EKF}) is found as follows:
\begin{equation}
\begin{array}{l}
\!\!s_{k+1}\!=\!s_k -\Delta\alpha_k s_k i_k + \Delta w_{sk}\\
\!\!i_{k+1}\!=\!i_k + \Delta\alpha_k s_k i_k - \Delta\beta i_k  + \Delta w_{ik}\\
\!\!\alpha_{k+1}\!=\!\alpha_k -\Delta\gamma \alpha_k + \Delta\gamma h[\mathbf{u}^*_k]+ \Delta w_{\alpha k}\\
\!\!\lambda_{1,k+1}\!=\!\lambda_{1k} + \Delta[\lambda_{1k}\!-\!\lambda_{2k}\!-1\!+\!\epsilon]\alpha_k i_k + \Delta\eta_{1k}\\
\!\!\lambda_{2,k+1}\!=\!\lambda_{2k} + \Delta[\lambda_{1k}\! -\!\lambda_{2k}\!-1\!+\!\epsilon]\alpha_k s_k + \Delta\beta\lambda_{2k} + \Delta\eta_{2k}\\
\!\!\lambda_{3,k+1}\!=\!\lambda_{3k} + \Delta[\lambda_{1k}\! -\!\lambda_{2k}\!-1\!+\!\epsilon]s_k i_k + \Delta\gamma\lambda_{3k} + \Delta\eta_{3k}\\
\!\!n_k=\alpha_k s_k i_k + v_k
\end{array}
\label{eq:EKFDiscrete}
\end{equation}
which can be formulated in a compact form:
\begin{equation}
    \begin{array}{l}
    \bm{s}_{k+1} = \mathbf{f}(\bm{s}_{k}, \bm{w}_k ; h(\mathbf{u}_k))\\
    n_k = g(\bm{s}_{k}) + v_k
    \end{array}
    \label{eq:statespacemodel}
\end{equation}
where $\mathbf{f}(\cdot)$ and $g(\cdot)$ represent the nonlinear equations in (\ref{eq:EKFDiscrete}). Following (\ref{eq:observationTotalCases}), if the number of confirmed cases is used as the observation, the second equation in (\ref{eq:statespacemodel}) is replaced by
\begin{equation}
c_k = s_0 - s_k + v_k = g(\bm{s}_{k}) + v_k
\label{eq:TotalCaseObservationsDiscrete}  
\end{equation}
which is a linear function of the state vector.

It is straightforward to linearize the discrete-time dynamic model (\ref{eq:EKFDiscrete}) by calculating its Jacobian matrices, as required for the implementation of the EKF/EKS. An alternative approach is to use a \textit{continuous-dynamics discrete-observations} approach, which is a classical method in optimal state estimation \cite{Gel74,AndersonMoore1979}. Accordingly, to implement the EKF/EKS, the state equations are updated by using the continuous version of the dynamic model (\ref{eq:EKF}), while the observations are only updated on discrete-time intervals (e.g., on a daily basis).

%\appendix[The Extended Kalman Filter Algorithm\label{sec:jacobians}]
%We define the discrete-time version of the dynamic model (\ref{eq:EKF}) variables as: $n_k = n(k\Delta)$, $s_k = s(k\Delta)$, $i_k = i(k\Delta)$ and $\alpha_k = \alpha(k\Delta)$, where $\Delta$ is the discretization time unit. Using a similar notation for the process and observation noises,
%The discretized version of the dynamic equations () can 
\section{The extended Kalman filter}
\label{sec:EKFEKS}
With the state vectors and discretized dynamic models defined in Appendix \ref{sec:discretized}, defining the noisy regular reports $x_k$ (number of new cases $n_k$, or alternatively the total confirmed cases $c_k$) with observation noise variance $r_k$, process noise covariance matrix $\mathbf{Q}_k$, initial state estimate $\hat{\bm{s}}_{0}^+$ with covariance $\mathbf{P}_{0}^+$, the recursive equations for the EKF are listed in Algorithm \ref{alg:EKF}. In this algorithm, $\mathbf{A}_k = \left.\partial \mathbf{f}/\partial \bm{s}_k\right|_{\bm{s} = \hat{\bm{s}}_{k}^+}$ and $\mathbf{c}_k = \left.\partial g/\partial \bm{s}_k\right|_{\bm{s} = \hat{\bm{s}}_{k}^-}$ are the linearized forms of the system's dynamics equations, and $\bar{\bm{w}}_k=\mathbb{E}\{\bm{w}_k\}$. Step\,10 of the algorithm corresponds to enforcing hard constraints on the estimated variables (e.g. positiveness or range) and Step\,11 corresponds to Kalman filter sanity checks (cf. Section~\ref{sec:performance_monitoring}).

\begin{algorithm}[tbh]
\begin{flushleft}
\caption{An extended Kalman filter for simultaneous compartment variable and model parameter tracking
\label{alg:EKF}}
\begin{algorithmic}[1]
\REQUIRE{The noisy regular reports $x_k$ (number of new cases $n_k$, or alternatively the total confirmed cases $c_k$)}
\REQUIRE Initial conditions: $\mathbf{Q}$, $\mathbf{R}$, $\hat{\bm{x}}_{0}^+$, $\mathbf{P}_{0}^+$
\ENSURE{$\hat{\bm{s}}_{k}^+$ (vector of state and model parameter estimates)}
\FOR{$k=0,\cdots, T$}
\STATE \textit{State prediction:}
\STATE $\hat{\bm{s}}_{k+1}^- = \mathbf{f}(\hat{\bm{s}}_{k}^+, \bar{\bm{w}}_k ; h(\mathbf{u}_k))$ 
\STATE $\mathbf{P}_{k+1}^- = \mathbf{A}_{k}\mathbf{P}_{k}^+\mathbf{A}_{k}^{T} + \mathbf{Q}_k$
\STATE \textit{Measurement update:}
\STATE $\mathbf{k}_{k} = \mathbf{P}_{k}^- \mathbf{c}_{k}[\mathbf{c}_{k}^T\mathbf{P}_{k}^-\mathbf{c}_{k} + r_k]^{-1}$
\STATE $\hat{v}_{k} = x_{k} - g(\hat{\bm{s}}_{k}^-)$
\STATE $\hat{\bm{s}}_{k}^+ = \hat{\bm{s}}_{k}^- + \mathbf{k}_{k}\hat{v}_{k}$
\STATE $\mathbf{P}_{k}^+ = [\mathbf{I} - \mathbf{k}_{k}\mathbf{c}_{k}^T]\mathbf{P}_{k}^-$
\STATE \textit{Check and enforce variable and parameter ranges using hard-constraints}
\STATE \textit{Performance monitoring}
\ENDFOR
\end{algorithmic}
\end{flushleft}
\end{algorithm}

\section*{Acknowledgment}
The author sincerely thanks Prof. Christian Jutten, Emeritus Professor of Universit\'e Grenoble Alpes (UGA) for the very insightful and motivating comments on the first versions of this work. The author also acknowledges the SEEPIA COVID-19 working group at UGA, directed by Prof. Didier Georges, for the fruitful (virtual) scientific meetings during the pandemic.  
%///////////////////////////////////////////////////////
%\newpage
\bibliographystyle{IEEEtran}
\bibliography{IEEEabrv,References}
%\newpage
%\listoffigures
%\listoftables
\end{document}